\documentclass[twocolumn]{autart_arxiv}    % Enable this line and disable the 
\usepackage{graphicx} 
\usepackage{subfigure} 
\usepackage{amsmath}        % Include this line if your 
\usepackage{amssymb}
\usepackage[dvips]{epsfig}    % or this line, depending on which
\usepackage{cases} 
\usepackage{color}
\usepackage{algorithm}
\usepackage{algorithmicx} 
\usepackage{pdfsync}
\usepackage{setspace} 
\usepackage{amscd,dsfont,texdraw,tikz,verbatim,color,multicol,multirow}
\usepackage[numbers]{natbib}

\usetikzlibrary{decorations.pathmorphing} % for snake lines
\usetikzlibrary{matrix} % for block alignment
\usetikzlibrary{calc} % for manimulation of coordinates
\usetikzlibrary{shapes,arrows,chains}                  % you prefer.

\newcommand{\bassume}{ \begin{assume} \begin{rm} }
\newcommand{\eassume}{ \end{rm} \end{assume} }
\newcommand{\bcondition}{ \begin{condition} \begin{rm} }
\newcommand{\econdition}{ \end{rm}  \end{condition} }
\newcommand{\bremark}{ \begin{remark} \begin{rm} }
\newcommand{\eremark}{ \end{rm} \hfill $\triangleleft$ \end{remark} }
\newcommand{\btheorem}{ \begin{theorem} \begin{rm} }
\newcommand{\etheorem}{ \end{rm} \hfill $\triangleleft$ \end{theorem} }
\newcommand{\blemma}{ \begin{lemma} \begin{rm} }
\newcommand{\elemma}{ \end{rm} \hfill $\triangleleft$ \end{lemma} }
\newcommand{\bcorollary}{ \begin{corollary} \begin{rm} }
\newcommand{\ecorollary}{ \end{rm}  \end{corollary} }
\newcommand{\bdefinition}{ \begin{definition}\begin{rm} }
\newcommand{\edefinition}{ \end{rm} \hfill $\triangleleft$ \end{definition} }
\newcommand{\bproposition}{ \begin{proposition} \begin{rm} }
\newcommand{\eproposition}{ \end{rm}  \end{proposition} }
\newcommand{\bexample}{ \begin{example} \begin{rm} }
\newcommand{\eexample}{ \end{rm} \hfill $\triangleleft$ \end{example} }
\newcommand{\bproblem}{ \begin{problem} \begin{rm} }
\newcommand{\eproblem}{ \end{rm} \hfill $\triangleleft$ \end{problem} }

\newcommand{\bproof}{ \textit{Proof:} \begin{rm} }
\newcommand{\eproof}{ \end{rm}  \hfill $\square$}

\newtheorem{theorem}{\bf Theorem}[section]
\newtheorem{lemma}{\bf Lemma}[section]
\newtheorem{definition}{\bf Definition}[section]
\newtheorem{remark}{\bf Remark}[section]
\newtheorem{corollary}{\bf Corollary}[section]
\newtheorem{proposition}{\bf Proposition}[section]
\newtheorem{example}{\bf Example}[section]
\newtheorem{assume}{\bf Assumption}[section]
\newtheorem{condition}{\bf Condition}[section]
\newtheorem{problem}{\bf Problem}
\begin{document}

\begin{frontmatter}
 \runtitle{Mean Field LQG Social Optimization: A Reinforcement Learning Approach}  % Running title for regular 

\title{Mean Field LQG Social Optimization: A Reinforcement Learning Approach}  %\thanksref{footnoteinfo} % Title, preferably not more 
% Model-Free Approach to Linear Quadratic Mean Field Games using Integral Reinforcement Learning                                    

% \thanks[footnoteinfo]{This paper was not presented at any IFAC 
% meeting. Corresponding author M.~T.~Cicero. Tel. +XXXIX-VI-mmmxxi. 
% Fax +XXXIX-VI-mmmxxv.}

\author[a]{Zhenhui Xu}\ead{zhenhuixu@sophia.ac.jp},
\author[b]{Bing-Chang Wang}\ead{bcwang@sdu.edu.cn}, and 
\author[a]{Tielong Shen}\ead{tetusin@sophia.ac.jp}
% \author[b]{Minyi Huang}\ead{mhuang@math.carleton.ca} % (ead) as shown

\address[a]{Department of Engineering and Applied Sciences, Sophia University, Tokyo 102-8554, Japan}
\address[b]{School of Control Science and Engineering, Shandong University, Jinan 250061, China}
          
                             % keyword list or with the 
\begin{abstract}
This paper presents a novel model-free method to solve linear quadratic Gaussian mean field social control problems in the presence of multiplicative noise. The objective is to achieve a social optimum by solving two algebraic Riccati equations (AREs) and determining a mean field (MF) state, both without requiring prior knowledge of individual system dynamics for all agents. In the proposed approach, we first employ integral reinforcement learning techniques to develop two model-free iterative equations that converge to solutions for the stochastic ARE and the induced indefinite ARE respectively. Then, the MF state is approximated, either through the Monte Carlo method with the obtained gain matrices or through the system identification with the measured data. Notably, a unified state and input samples collected from a single agent are used in both iterations and identification procedure, making the method more computationally efficient and scalable. Finally, a numerical example is given to demonstrate the effectiveness of the proposed algorithm.
\end{abstract}
\end{frontmatter}

\section{Introduction}
The linear quadratic Gaussian social mean field control provides a powerful tool for solving cooperative differential games, which are characterized by a large population of agents operating under linear dynamics and seeking to optimize a common quadratic cost function. Its analytical tractability has attracted considerable attention and found application in various practical fields, such as opinion dynamics \cite{wang2022robust}, production adjustment \cite{wang2019mean}, and dynamic economic models \cite{wang2023rational}, etc. 
 % a large population of agents, each following linear dynamics, interacts in a way that overall performance.}
\subsection{MFG theory and social optimum} 
As a branch of game theory, mean field game (MFG) theory is devoted to the analysis of decentralized decision problems involving a large number of interacting agents. In this area, the seminal contributions of Lasry and Lions \cite{lasry2006jeux,lasry2007mean} and Huang {\sl et al.} \cite{huang2006large,huang2007large}, initially influential in theoretical research, have not only seen extensions but also found widespread practical applications (see, e.g. \cite{chen2015state,bauso2016opinion,siwe2016network,kuga2019vaccinate}). The effectiveness of MFG  in solving large-scale system problems lies in its basic idea of decoupling the interaction between agents by using a mean field representation. As a result, solutions to MFGs can be characterized by a system of forward-backward equations, where the forward Kolmogorov-Fokker-Planck equation represents the evolution of the mean field and the backward Hamilton-Jacobi-Bellman equation describes the evolution of the value function of an optimal control problem. For further analysis details of MFGs, the interested readers are referred to \cite{bensoussan2013mean} for an overview of the MFGs and mean field type control, \cite{cousin2011mean} for an introduction of both theory and applications of MFGs, \cite{carmona2018probabilistic} for a probabilistic analysis of MFGs, and \cite{gomes2014mean} for a comprehensive survey of MFGs.
 
 Besides non-cooperative games, social optima in MF models has received significant attention. Here, $N$ cooperative agents collaborate to minimize a social cost, which is the sum of all agents' costs. This problem can be viewed as a team decision
problem (\cite{marschak1955elements}) that has been studied extensively over time. The MF control methodology offers a key advantage in this kind of problems, allowing us to analyze the impact of an individual agent's strategy on the social cost while keeping the strategies of other agents fixed. One approach, as described in \cite{huang2012social}, transforms a social optimization problem into a non-cooperative game via a person to person optimization approach. It leads to a set of decentralized person-by-person optimal strategies, which serves as an approximation to the social optimum. Based on this idea, subsequent studies have explored various aspects of MF models. For example, \cite{huang2016linear} focused on a linear quadratic Gaussian (LQG) model involving a major agent and a large number of minor agents, and \cite{wang2017social} studied dynamic team decision problems where a Markov jump parameter appears as a common source of randomness.  For more literature on the topic, see \cite{arabneydi2015team,salhab2018dynamic,du2022social}. In addition, solving MF social optimization problems has been approached through the direct approach \cite{huang2019linear2,huang2021linear2}. The direct approach involves solving a game among $N$ players directly, resulting in a system of coupled equations, which is then analyzed in the limit as the number of agents $N$ approaches infinity. This approach has also yielded fruitful results, especially in the LQG setting, leading to closed-form solutions and theoretical tractability. Such works include \cite{wang2020mean}, which considered the LQG control systems with indefinite state weight matrix, \cite{wang2020indefinite}, which involved the multiplicative noise, and \cite{wang2022linear}, which incorporated the common noise.

\subsection{Model-based and model-free numerical methods}
The explicit resolution of MFG problems is generally difficult. Consequently, researchers have directed their efforts towards the development of numerical methods to tackle various types of MF models. At the beginning, numerical methods were designed for low-dimensional MFGs by using finite differences \cite{achdou2010mean}, semi-Lagrangian schemes \cite{carlini2014fully}, and primal-dual methods \cite{briceno2018proximal}. More advanced techniques, such as neural network approximations \cite{ruthotto2020machine,fouque2020deep,carmona2021convergence} and fixed-point iterations \cite{cacace2021policy,camilli2022rates,lauriere2023policy}, have been developed to solve high-dimensional MFGs. However, these results relied on exact knowledge of agents' dynamics, limiting their practical applicability.

 Recently, reinforcement learning (RL) has emerged as a very promising approach for solving MFG problems in a model-free manner. These methods involve alternating updates between the MF state and the control policy. Typically, in each iteration, the policy is updated by computing a best response against the mean field. Various techniques have been developed to implement this idea, including fictitious play \cite{elie2020convergence,elie2019approximate}, Q-learning \cite{guo2019learning,anahtarci2023learning}, and deep reinforcement learning \cite{cui2021approximately,perrin2021mean}. Alternatively, evaluating a policy directly instead of searching for an optimal policy is another useful approach \cite{lauriere2022scalable}. To address the challenge of updating the MF state in scenarios involving large finite or continuous state space, researchers have introduced neural network approximations to represent the mean field evolution (\cite{perrin2021mean}). Additionally, Monte Carlo samples have been used to implement an empirical distribution (\cite{angiuli2022unified}). We stress that our main focus lies on the LQG setting, in which case the MF state dynamics exhibit a closed-loop linear structure. \cite{uz2020reinforcement} exploited this structural property of the MF update operator in the LQG case and devised a double-loop method using an actor-critic algorithm.  \cite{uz2023reinforcement} solved non-stationary LQG-MFGs involving multiple populations by employing zero-order stochastic optimization. Furthermore, \cite{xu2023mean} computed a set of decentralized strategies for continuous-time LQG-MFGs by using trajectories of a single agent.

\subsection{Motivation and contribution}
Most of the above literature focuses on the application of model-free methods to non-cooperative mean-field games, while there is a limited research on model-free methods for MF social optimization problems. Relevant results can be found in \cite{gu2021mean,carmona2019linear,angiuli2022unified}. In \cite{gu2021mean}, a model-free kernel-based Q-learning algorithm with a linear convergence rate was proposed. \cite{angiuli2022unified} further replaced a MF simulator with estimation and developed a Q-learning method. \cite{carmona2019linear} studied this problem in the discrete time setting using the McKean-Vlasov simulator and policy gradient method. All these results considered additive noise. However, in practical situations, multiplicative noise is more realistic description of stochastic disturbances as the intensity noise is dependent on the state and input \cite{yong1999stochastic,wang2019mean1,wang2020indefinite,moon2020linear,aberkane2023addendum}. Therefore, there remains a clear need for further investigation into model-free methods for the mean field social optimization problems with mutiplicative noise, which motivates the present study.

In this paper, we study the continuous-time LQG-MF social optimization problems with multiplicative noise. If prior information on the system dynamics is available, the optimal gains can be determined by solving two algebraic Riccati equations, while the corresponding mean field state is obtained through an ordinary differential equation (\cite{wang2020indefinite}). However, the technical challenge arises when the dynamics of agents are completely unknown. First, the existence of  multiplicative noise, ({\sl i.e., the diffusion term in the agents' dynamics depends on state and control variables}), results in a stochastic ARE (SARE). This equation cannot be solved using existing model-free methods. Second, the other generalized ARE involves a term composed of the solution of the SARE and the diffusion coefficient matrices, which hinders the use of traditional reinforcement learning methods (\cite{lee2012integral,jiang2012computational}). Therefore, solving this equation calls for a new model-free algorithm. 
 
To overcome these hurdles, a novel model-free design scheme has been developed. This scheme adopts two off-policy RL methods to sequentially obtain the solutions of two AREs, and uses two approximate methods to calculate the MF state corresponding to the MF social optimum. In the first step, an integral reinforcement learning technique is utilized to construct a model-free iterative equation that approximates the SARE. In contrast to \cite{li2022stochastic}, where the proposed method only removes the need for the state matrix and relies on the state weight matrix as a positive definite matrix, our method goes to a step further by completely eliminating the requirement for all system parameters from the iterative equation and relaxing the state weight matrix to a positive semidefinite matrix. The second step focuses on solving the other ARE. An equation transformation is initially applied to remove the system dynamics from the weight matrix of the ARE, resulting in an indefinite ARE. To address this deduced equation, a model-based policy iteration is designed with the guaranteed convergence. Furthermore, an iterative equation is developed, completely removing the requirement of system matrices. Finally, based on the obtained gain matrices, the MF state can be approximated either by using Monte Carlo samples from a single agent or by identifying the system parameters. By establishing the equivalence between the model-free methods and the model-based methods under the achievable rank condition, we ensure the convergence of the RL algorithms.
% The paper provides rigorous proofs that all off-policy iterations are convergent under achievable conditions.

The novelty of the paper can be outlined as follows:
\begin{enumerate}
    \item[1).] We develop a model-free scheme that does not require any system matrices and can efficiently determine the solutions of the SARE and the indefinite ARE. Especially, for solving the indefinite ARE, the paper establishes the Hurwitz property of the closed-loop matrix for iterative equations involving an indefinite state weight matrix, which differs from the traditional policy iteration for linear quadratic optimal control problems.
    \item[2).] We introduce a system identification-based MF approximation method. In contrast to the Monte Carlo method, data gathered from RL enables the estimation of partial parameters for individual dynamical models, eliminating the necessity for data re-collection in system identification-based MF calculations. 
    \item[3).] For the model-free decentralized strategy design, both policy iterations and system identification procedures are performed using a unified dataset sampled from the agent’ s state trajectories and inputs, enhancing computational efficiency. 
\end{enumerate}
 % and demonstrates that all policy iterations can be carried out utilizing a unified dataset sampled from the agent's state trajectories and inputs, enhancing computational efficiency.

After providing a list of notations, the rest parts of this paper are structured as follows. Section \ref{sec2} describes the problem and presents basic results for LQG mean field social optimization. Section \ref{sec3} presents a model-free approach to compute the mean field social optimal decentralized strategy using real-time samples from a single agent. In Section \ref{sec4}, a simulation example is presented to demonstrate the effectiveness of the proposed approach. Finally, Section \ref{sec5} concludes the paper.

\subsection{Notations}
  A list of notations is presented as follows. $\mathbb{N}_+$: the set of all positive integers; for a family of $\mathbb{R}^n$-valued random variables $\{x(\tau),\tau\geq0\}$, $\sigma(x(\tau),\tau\leq t)$ is the $\sigma$-algebra generated by these random variables; ${\bf 0}$: the set of all zero matrices (vectors) having compatible dimensions; for $v\in\mathbb{R}^m$ and $P\in\mathbb{R}^{m\times m}$, $\|v\|_P^2=v^{\mathrm{T}}Pv$; $\otimes$: the Kronecker product; $e_j$: an $n$-dimensional row vector with the $j$-th element being one and the rest being zero; $E_j$: an $n\times n$-dimensional matrix with the $j$-th diagonal element being one and the rest being zero; $\mathbb{S}^n$: the set of all $n\times n$ symmetric matrices; for $R\in\mathbb{S}^n$, $R>0$: $R$ is a positive definite matrix; for $A\in\mathbb{S}^n$ and $B\in\mathbb{S}^n$, $A>B$($A<B$): $A-B>0$($A-B<0$); $\mathrm{col}(A)$: the $m n$-dimensional vector formed by stacking the columns of $A\in\mathbb{R}^{m\times n}$ on top of one another; $L_2([0,\infty),\mathbb{R}^n)\!=\!\{f\!:\![0,\infty)\!\rightarrow\!\mathbb{R}^n| \int_0^{\infty}\!\|f(\tau)\|^2\mathrm{d}\tau\!\!\!<\!\!\infty\}$; for $P\!\in\!\mathbb{S}^{n}$ and $x\!\in\!\mathbb{R}^n$, define the following operators:  $\bar{P}=[p_{11},2p_{12},\cdots,2p_{1n},p_{22},2p_{23},\cdots,2p_{n-1,n},p_{nn}]^{\mathrm{T}}$ $\in$$\mathbb{R}^{\frac{1}{2}n(n+1)}$ and $\hat{x}=[x_1^2,x_1x_2,\cdots,x_1x_n,x_2^2,\cdots,x_{n-1}x_n,$ $x_n^2]^{\mathrm{T}}$$\in$$\mathbb{R}^{\frac{1}{2}n(n+1)}$.

\section{Problem Description and Preliminaries}\label{sec2}
 We consider a population of $N$ agents, denoted as $\mathcal{A}=\{\mathcal{A}_i,1\leq i\leq N\}$, where $\mathcal{A}_i$ represents the $i$-th agent. The state process $x_i(t)$ for $\mathcal{A}_i$ satisfies the following stochastic differential equation (SDE)
\begin{equation}\label{sys1}
\mathrm{d}x_i(t)  \!= \!(\!Ax_i(t)+Bu_i(t)\!)\mathrm{d}t\!+\!(\!Cx_i(t)\!+\!Du_i(t)\!)\mathrm{d}w_i(t),
\end{equation}
where the vectors $x_i\in\mathbb{R}^n$ and $u_i\in\mathbb{R}^m$ represent the state and the input, respectively. The noise processes $\{w_i,1\leq i\leq N\}$ are $N$ independent standard one-dimensional Brownian motions defined on a compact filtered probability space $(\Omega,\mathcal{F},\{\mathcal{F}_t\}_{0\leq t\leq T},\mathbb{P})$. The constant real matrices $A,B,C$, and $D$ all have compatible dimensions. Denote the control input of the population $\mathcal{A}$ as $u=\{u_1,\cdots,u_N\}$.

The social optimization problem consists of $N$ agents' dynamics (\ref{sys1}) and the following social cost
\vspace{-0.1cm}
\begin{equation}\label{Jsoc}
J_{\mathrm{soc}}(u)=\sum_{i=1}^{N}J_i(u),
\vspace{-0.3cm}
\end{equation}
\vspace{-0.1cm}
where 
\vspace{-0.1cm}
\begin{align*}
J_i\big(u\big)\!=\!E \left[\int_0^{\infty}\!\!\big(\|x_i- \Gamma x_{(N)}\|_{Q}^{2}\!+\!\|u_i\|_{R}^2 \big) \mathrm{d}\tau\!\right],\label{J}
\end{align*}
 $x_{(N)}=\frac{1}{N}\sum_{j=1}^Nx_j$ is the average state of $\mathcal{A}$, the weight matrices $Q\in\mathbb{S}^n$ and $R\in\mathbb{S}^{m}$ satisfy $Q\geq0$ and $R>{0}$, and $\Gamma\in\mathbb{R}^{n\times n}$ . 

 Throughout the paper, It is aimed to find a social optimum for (\ref{sys1})-(\ref{Jsoc}). The related concepts of stability, observability, admissibility, and optimality are introduced as follows. 

We first consider the following system:
 \begin{numcases}
\!\!\mathrm{d}y(t) \!=\! \!(\!Ay(t)\!\!+\!\!B\!u(t)\!)\mathrm{d}t\!+\!\!(\!Cy(t)\!\!+\!\!Du(t)\!)\mathrm{d}w(t),t\geq0, \label{sysdef1}\\
\!\!z(t) = Fy(t),\label{sysdef2}
 \end{numcases}
where the initial state $y(0)=y_0\in\mathbb{R}^n$, and $w(t)$ is a one-dimensional Brownian motion. 
 \begin{definition}\label{A0}
The system (\ref{sysdef1}) (or simply $[A,B;C,D]$) is said to be mean square stabilizable, if there exists a constant matrix $K\in\mathbb{R}^{m\times n}$ such that the (unique) strong solution of 
\begin{equation*}\!\mathrm{d}y(t) \!= \!(\!A\!-\!BK)y(t)\mathrm{d}t\!+\!(\!C\!-\!DK)y(t)\mathrm{d}w(t)
\end{equation*}
satisfies $\lim_{t\rightarrow\infty}E[\|y(t)\|^2]=0$ for all $y_0\in\mathbb{R}^n$. In this case, $K$ is called a stabilizer of $[A,B;C,D]$.
\end{definition}
\begin{definition}[\cite{wang2020indefinite}]
The system (\ref{sysdef1})-(\ref{sysdef2}) ( or simply $[A,C;F]$) is said to be exactly observable, if there exists a $T_0\geq0$ such that for any $T\geq T_0$, $z(t)={\bf 0}$, $u(t)={\bf 0}$, {\sl a.s., $ 0\leq t\leq T$ implies $y_0={\bf 0}$}. If the diffusion term of the stochastic system (\ref{sysdef1}) becomes $\sum_{j=1}^k(C_jx+D_ju)\mathrm{d}w_i$, then the system, denoted by $[A,C_1,C_2,\cdots,C_k;F]$, is exactly observable. If $C={\bf0}$, then the system, abbreviated as $(A,F)$, is observable.
\end{definition}
\begin{definition} A decentralized control policy of $\mathcal{A}_i$ is said to be admissible, denoted by $u_i\in\mathcal{U}_{di},$ if
\begin{enumerate}
\item[1).] $u_i$ is adapted to $\sigma(x_i(\tau),0\leq\tau\leq t)$; \\
\item[2).] $E[\int_{t}^{\infty}\|x_i(\tau)\|^2\mathrm{d}\tau]<\infty,~~\forall t\geq0$.
\end{enumerate}
\end{definition}
\begin{definition} A centralized control policy of $\mathcal{A}_i$ is said to be admissible, denoted by $u_i\in\mathcal{U}_{ci},$ if 
\begin{enumerate}
\item[1).] $u_i$ is adapted to $\sigma(\bigcup _{i=1}^N \mathcal{F}_t^i)$, $\mathcal{F}_t^i=\sigma(x_i(0),w_i(\tau),0\leq \tau\leq t)$, $i=1,\cdots,N$;\\
\item[2).] $E[\int_{t}^{\infty}\|x_i(\tau)\|^2\mathrm{d}\tau]<\infty,~~\forall t\geq0$.
\end{enumerate}
\end{definition}
 \begin{definition}[\cite{wang2020social,wang2020indefinite}]\label{def2.4} A sequence of decentralized control policies ${u}^o\!=\{{u}^o_1,\cdots,{u}^o_N\}$ for $\mathcal{A}$ is said to have asymptotic social optimality if 
\begin{equation*}
\left|\frac{1}{N}J_{\mathrm{soc}}({u}^o)-\frac{1}{N}\inf_{u\in\cup_{i=1}^N\mathcal{U}_{ci}}J_{\mathrm{soc}}(u)\right|=o\left(1\right).
\end{equation*}
\end{definition}
We define $\Upsilon\triangleq R+D^{\mathrm{T}}PD$ and $Q_{\Gamma}\triangleq \Gamma^{\mathrm{T}}Q+Q\Gamma-\Gamma^{\mathrm{T}}Q\Gamma$, where $P\in\mathbb{S}^{n}$ is a solution of the following ARE:
\begin{equation}\label{ARE1}
\begin{aligned}
&A^{\mathrm{T}}P+PA+C^{\mathrm{T}}PC-\left(B^{\mathrm{T}}P+D^{\mathrm{T}}PC\right)^{\mathrm{T}}\Upsilon^{-1}\\
&\times\left(B^{\mathrm{T}}P+D^{\mathrm{T}}PC\right)+Q ={\bf0}.
\end{aligned}
\end{equation}
Let $\Pi\in\mathbb{S}^n$ be a solution of 
\begin{equation}\label{ARE2}
\begin{aligned}
&A^{\mathrm{T}}\Pi+\Pi A\!-\!\left(B^{\!\mathrm{T}}\Pi\!+\!D^{\!\mathrm{T}}PC\right)^{\mathrm{T}}\Upsilon^{-1}\!(B^{\!\mathrm{T}}\Pi\!+\!D^{\!\mathrm{T}}PC)\\
&+C^{\mathrm{T}}PC+Q-Q_{\Gamma}={\bf0},
\end{aligned}
\end{equation}
and let $\bar{x}\in L_2([0,\infty),\mathbb{R}^n)$ be a mean field state satisfying
\begin{equation}\label{xbar}
\!\!\!\dot{\bar{x}}(t)\!=\!\left(\!A\!-\!B\Upsilon^{-1}\left(B^{\mathrm{T}}\Pi+D^{\mathrm{T}}PC\right)\right)\bar{x}(t),~\bar{x}(0)\!=\!\bar{x}_{0}.
\end{equation} 
According to \cite[Theorem 4.4]{wang2020indefinite}, it is easy to show that the set of decentralized control polices given by
\begin{equation}\label{uhat}
\begin{aligned}
{u}^o_i=&-\Upsilon^{-1}(B^{\!\mathrm{T}}P+D^{\mathrm{T}}PC)x_i \\
&-\Upsilon^{-1}B^{\mathrm{T}}(\Pi-P)\bar{x},~~i=1,\cdots,N,
\end{aligned}
\end{equation}
has asymptotic social optimality under the following basic assumptions:
% \vspace{-0.5cm}
\begin{enumerate}
\item[({\bf{A1}})] $x_i(0)$, $i=1,\cdots,N$ are mutually independent and also independent of $\{w_i,i=1,\cdots,N\}$. They share the same mathematical expectation, {\sl i.e., $Ex_i(0)=\bar{x}_0$,} and have a finite second moment {\sl i.e., $\max_{1\leq i\leq N}E\|x_i(0)\|^2<c_0$, where $c_0<\infty$ is a constant independent of $N$}.
\item[({\bf{A2}})] The system $[A,B;C,D]$ is mean square stabilizable.
\item[({\bf{A3}})] $[A,C;\sqrt{Q}]$ is exactly observable, and $(A,\sqrt{Q}(I-\Gamma))$ is observable. 
\end{enumerate}  
\begin{remark}\label{rem2.1}
Assumptions ({\bf A2}) and ({\bf A3}) are standard in the stochastic control literature (see, e.g. \cite{li2022stochastic,wang2020indefinite,zhang2004stabilizability}). Specifically, mean square stability and exact observability (an extended version of observability) are crucial for solving stochastic control problems involving multiplicative noise. Under assumptions (A1) and (A2), there exists a unique solution $P>0$ to equation (\ref{ARE1}) such that $\Upsilon^{-1}(B^{\mathrm{T}}P+D^{\mathrm{T}}PC)$ is a stabilizer of $[A,B;C,D]$. Moreover, equation (\ref{ARE2}) admits a unique positive definite solution $\Pi$ .
\end{remark}
Furthermore, as stated in Definition \ref{def2.4}, the optimal centralized solution serves as a benchmark for comparison against the decentralized control policies (\ref{uhat}). To ensure the existence of an optimal centralized solution, we need an additional assumption:

({\bf{A4}}) $[A,{\bf C}_1,{\bf C}_2,\cdots,{\bf C}_N;\sqrt{{\bf \bar{Q}}}]$ is exactly observable, where ${\bf C}_j={E}^N_j\otimes C$, $j=1,\cdots,N$, and ${\bf \bar{Q}}\!=\!(\bar{Q}_{ij})$, $\bar{Q}_{ii}\!=\!Q\!-\!Q_{\Gamma}/N$, $\bar{Q}_{ij}\!=\!-Q_{\Gamma}/N$, $1\leq i\neq j\leq N$.

\begin{remark}
In \cite{wang2020indefinite,zhang2004stabilizability}, the derivation of optimal centralized solution involves employing augmentation state space technique to get a compact representation of the dynamics of all agents 
\begin{equation}\label{sysaug}
\begin{aligned}
\mathrm{d}{\bf x}=\left({\bf\check{A}}{\bf x}\!+\!{\bf B}{\bf u}\right)\mathrm{d}t+\sum_{j=1}^{N}\left({\bf C}_j{\bf x}+{\bf D}_j{\bf u}\right)\mathrm{d}w_j,
\end{aligned}
\end{equation}
and solving the coupled SARE
\begin{equation}\label{Jaug}
\begin{aligned}
&\!\!\!\!\!\!\!{\bf{\check{A}}}^{\!\!\mathrm{T}}{\bf P}\!+\!{\bf P}{\bf{\check{A}}}\!+\!\sum_{i=1}^{N}{\bf C}_i^{\!\mathrm{T}}{\bf P C}_i\!+\!{\bf \bar{Q}}\!-\!\!\left(\!\!{\bf B}^{\mathrm{T}}{\bf P}\!+\!\sum_{i=1}^N{\bf D}_i^{\!\mathrm{T}}{\bf P}{\bf C}_i\!\!\right)^{\!\!\mathrm{T}}\\
&\!\!\!\!\!\!\!\times {\bf \Upsilon}^{-1}\!\!\left({\bf B}^{\mathrm{T}}{\bf P}+\sum_{i=1}^{N}{\bf D}_i^{\mathrm{T}}{\bf PC}_i\right)={\bf0},
\end{aligned}
\end{equation}
where ${\bf x}\!=\![x_1^{\mathrm{T}},\!\cdots\!,x_N^{\mathrm{T}}]^{\mathrm{T}}$, ${\bf u}\!=\![u_1^{\mathrm{T}},\!\cdots\!,u_N^{\mathrm{T}}]^{\mathrm{T}}$, ${\bf \check{A}}\!=\!\mathrm{diag}(A,\cdots,A)$, ${\bf B}\!=\!\mathrm{diag}(B,\cdots,B)$, ${\bf D}_j={E}^N_j\otimes D$, $j=1,\cdots,N$, and $\!\cdots\!,{\bf0})$, ${\bf R}\!=\!\mathrm{diag}(R,\!\cdots\!,R)$.

The existence of a stabilizing solution ${\bf P}>0$ for SARE (\ref{Jaug}) is guaranteed by ({\bf A4}) together with the following facts: 1). the mean square stability of system (\ref{sysaug}) directly follows from that of $[A,B;C,D]$; 2). $Q\geq0$ implies the positive semidefiniteness of $\bar{\bf{Q}}$. Additionally, in the case of $Q>0$, ({\bf A4}) can be simplified to $\Gamma\neq \frac{1}{N}I_n$. As a result, the social optimization problem admits a centralized solution (\cite{wang2020indefinite,zhang2004stabilizability}). 
\end{remark} 

The meaning of the result summarized above is that the decentralized social optimization problem might be found by solving AREs (\ref{ARE1})-(\ref{ARE2}) and ODE (\ref{xbar}) when the dynamical model (\ref{sys1}) are known in advance. However, the requirement on exact knowledge of agent's dynamics confronted by a lot of theoretical optimization algorithm is a bottleneck in the practice. Consequently, the primary objective of this paper is to address this limitation by introducing a model-free algorithm capable of solving the decentralized social optimization problem. 

% The main purpose of this paper is to develop a model-free algorithm for solving the decentralized social optimization problem.
% The meaning of the result summarized above is that the decentralized social optimization problem might be solved by coupled AREs (\ref{ARE1})-(\ref{ARE2}) and the dynamics of MF state (\ref{xbar}), if the system dynamics (\ref{sys1}) are exactly known previously. However, the prior condition on exact knowledge of agent's dynamical model confronted by a lot of theoretical optimization algorithm is a bottleneck in the practice. The main purpose of this paper is to develop a model-free algorithm for solving the decentralized social optimization problem. 
% \linespread{0.5}
\section{Model-free Decentralized Control Design}\label{sec3}
In this section, we establish a new model-free method to approximate the set of decentralized control policies (\ref{uhat}). First, two model-based policy iteration (PI) algorithms are proposed to approximately solve the two AREs. Next, by employing the RL technique, we remove their dependence on system dynamics. Finally, based on the obtained gain matrices, we calculate a mean field state approximation.
\subsection{Model-based policy iteration}
For notational simplicity, let
\begin{equation}\label{kop}
K=\Upsilon^{-1}(B^{\mathrm{T}}P+D^{\mathrm{T}}PC).
\end{equation}
% Consider the condition ({\bf A2}) and $Q>0$, according to \cite{rami2000linear}, it is easy to show that there exists a unique positive definite solution $P$ for the SARE (\ref{ARE1}), and $K$ is a stabilizer for system $[A,B;C,D]$. However, solving the SARE (\ref{ARE1}) analytically requires system information and is \textcolor{blue}{often} difficult. 
Denote the $k$-th iteration of $K$ as
\begin{equation}\label{kk}
\!K_{k}\!=\!(R\!+\!D^{\mathrm{T}}P_{k}D)^{\!-1}(B^{\!\mathrm{T}}\!P_{k}\!+\!D^{\!\mathrm{T}}P_{k}C),~~k\in\mathbb{N}_+ 
\end{equation} 
with the $k$-th iterative solution $P_k$ of
\begin{equation}\label{pk}
\begin{aligned}
A_{k-1}^{\mathrm{T}}P_{k}+P_{k}A_{k-1}+C_{k-1}^{\mathrm{T}}P_kC_{k-1}+Q_{k-1}={\bf0},
\end{aligned}
\end{equation}
correspondingly, $A_{k}\triangleq A-BK_{k}$, $C_{k}\triangleq C-DK_{k}$, and $Q_{k}\triangleq Q+K_{k}^{\mathrm{T}}RK_{k}$.

The following lemma provides a PI method to numerically approximate $P$ and $K$.
\begin{lem}\label{lem1}
Suppose there exists a stabilizer $K_0\in\mathbb{R}^{m\times n}$ of $[A,B;C,D]$. Let $P_k\in\mathbb{S}^n$ be a solution of (\ref{pk}) and  $K_{k}$ be recursively updated by (\ref{kk}), then $\lim_{k\rightarrow\infty}P_{k}=P$ and $\lim_{k\rightarrow\infty}K_{k}=K$.
\end{lem}
See the proof in Appendix \ref{pf0}.
% \begin{remark}\label{rem3.1}
% As mentioned in the proof of Lemma \ref{lem1}, $(P,K)$ satisfies (\ref{pk}), which implies that equation (\ref{ARE1}) can be rewritten as the following Lyapunov function
% \begin{equation*}
% (A-BK)^{\mathrm{T}}P+P(A-BK)=-\mathbb{Q},
% \end{equation*}
% where $\mathbb{Q}=Q+K^{\mathrm{T}}RK+(C-DK)^{\mathrm{T}}P(C-DK)$. Since $P,Q,R>0$, the closed-loop system matrix $A-BK$ has all eigenvalues with negative real parts. This conclusion will be used in the following derivation.
% \end{remark}
% Given the conditions that $(P,K)$ satisfies (\ref{pk}) and $P,Q>{0}$, it can be shown that $A-BK$ is Hurwitz. Furthermore, by considering the facts $B\Upsilon^{-1}B^{\mathrm{T}}\geq{0}$, $\Upsilon>{0}$, and $Q_{\Gamma}=Q_{\Gamma}^{\mathrm{T}}$, and applying \cite[Theorem 2.1]{gohberg1986hermitian} along with Assumption \ref{As}, it can be concluded that there exists a unique maximal solution $S$ such that all eigenvalues of $A-BK-B\Upsilon^{-1}B^{\mathrm{T}}S$ lie in the closed left halfplane.

Before developing a second policy iteration algorithm, we define
\begin{equation}\label{trans1}
S\triangleq\Pi-P. 
\end{equation}
Combining equations (\ref{ARE1}) and (\ref{ARE2}) yields that $S$ satisfies
\begin{equation}\label{ARE2-2}
(A-BK)^{\!\mathrm{T}}S\!+\!S(A-BK)\!-\!SB\Upsilon^{\!-1}B^{\!\mathrm{T}}S\!-\!Q_{\Gamma}\!=\!{\bf0}.
\end{equation}
Since equations (\ref{ARE1}) and (\ref{ARE2}) have solutions, the symmetric solution set of (\ref{ARE2-2}) is not empty.
% \begin{remark}
% \textcolor{blue}{The Hurwitz property of the matrix $A-BK$ follows from the fact that $K$ is a stabilizer of $[A,B;C,D]$.} Furthermore, according to $B\Upsilon^{-1}B^{\mathrm{T}}\geq{0}$, $\Upsilon>{0}$, and $Q_{\Gamma}=Q_{\Gamma}^{\mathrm{T}}$, and using \cite[Theorem 2.1]{gohberg1986hermitian}, it can be deduced that a unique maximal solution $S$ exists such that all eigenvalues of $A-BK-B\Upsilon^{-1}B^{\mathrm{T}}S$ lie in the closed left halfplane.
% \end{remark}
\begin{remark}\label{rem3.3}
Compared to (\ref{ARE2}), equation (\ref{ARE2-2}) exhibits more pronounced features, leading to the following conclusions that facilitate the development of the subsequent iterative equations:
\begin{enumerate}
\item[1).] The Hurwitz property of matrix $A-BK$ arises from the fact that $K$ is a stabilizer of $[A,B;C,D]$. Therefore, the initial condition $K_s^0$ can be easily set to  ${\bf0}$.
\item[2).] The terms $C^{\!\mathrm{T}}\!P C$ and $D^{\!\mathrm{T}}\!P C$ disappear in (\ref{ARE2-2}), eliminating the use of the system information ($C$ and $D$) in the calculation of $S$. 
\end{enumerate}
Therefore, our objective now shifts from approximating the maximal solution of equation (\ref{ARE2}) to finding a sequence that converges to the maximal solution of equation (\ref{ARE2-2}).
\end{remark}

% Therefore, our next goal is not approximating the solution of equation (\ref{ARE2}), but to find a sequence that converges to the solution of (\ref{ARE2-2}).
To proceed, let 
\begin{equation}\label{ks}
K_s = \Upsilon^{-1}B^{\mathrm{T}}S,
\end{equation}
and define $\tilde{A}_{k}\triangleq A-B(K+K_s^{k})$, $Q_{\Gamma}^{k}\triangleq -Q_{\Gamma}+(K_s^{k})^{\mathrm{T}}\Upsilon K_s^{k}$. Then, we present a policy iteration algorithm in the following lemma.
\begin{lem}\label{lem3.2}
Let $S_k\in\mathbb{S}^n$ be the solution of the following Lyapunov equation
\begin{equation}\label{sk}
\tilde{A}_{k-1}^{\mathrm{T}}S_k+S_k\tilde{A}_{k-1}+{Q}^{k-1}_{\Gamma}={\bf0},
\end{equation}
and $K_s^{k}$ be recursively updated by
\begin{equation}\label{ksk}
K_s^k=\Upsilon^{-1}B^{\mathrm{T}}S_k,~~k\in\mathbb{N}_+.
\end{equation}
If $K_s^0={\bf0}$, then the sequence $\{S_k,K_s^k\}_1^{\infty}$ satisfies the following properties for all $k\in\mathbb{N}_+$.
\vspace{-0.1cm}
\begin{enumerate}
\item[1).] $S\leq S_{k+1}\leq S_{k}$. 
\item[2).] $\tilde{A}_{k-1}$ is Hurwitz.
\item[3).] $\lim_{k\rightarrow\infty}S_k=S$ and $\lim_{k\rightarrow\infty}K_s^k=K_s$.
\end{enumerate}
\end{lem}
\begin{pf}
We first show that the inequality $S_k\geq S$ and the result 2) hold for all $k\in\mathbb{N}_+$ by mathematical induction.
\vspace{-0.5cm}
\begin{enumerate}
\item[\romannumeral1).] If $k=1$, by the fact that $A-BK$ is Hurwitz, the initial condition $K_s^0={\bf0}$ ensures that $\tilde{A}_0$ is Hurwitz.

To show the inequality holds, we rewrite equation (\ref{ARE2-2}) as
\begin{equation}\label{lem2pfeq1}
\begin{aligned}
&(A-BK)^{\mathrm{T}}S+S(A-BK)-K_s^{\mathrm{T}}\Upsilon K_s-Q_{\Gamma}={\bf0},
\end{aligned}
\end{equation}
and subtract equation (\ref{lem2pfeq1}) from equation (\ref{sk}) with $k=1$, yielding
\begin{equation*} 
\begin{aligned}
\!\!\!(S_1-S)\tilde{A}_0+\tilde{A}_0^{\mathrm{T}}(S_1-S)=-(K_s)^{\mathrm{T}} \Upsilon K_s.
\end{aligned}
\end{equation*}
It follows from the Hurwitz property of $\tilde{A}_0$ and \cite[Proposition 8.13.1]{lancaster1985theory} that 
\begin{equation*}
\begin{aligned}
S_1\!-\!S\!=\!\int_0^{\infty}\!\!\!e^{\tilde{A}_0\tau}(K_s)^{\mathrm{T}}\Upsilon K_s e^{\tilde{A}_0^{\mathrm{T}}\tau}\mathrm{d}\tau,
\end{aligned}
\end{equation*}
which together with $\Upsilon>0$ implies $S_1\geq S$.
 
\item[\romannumeral2).] Suppose the result 2) holds for $k>1$, we prove that the result 2) also holds for $k+1$. Furthermore, the inequality $S_k\geq S$ holds for $k>1$. 

To this end, we subtract and add $(K_s^{k-1})^{\mathrm{T}}B^{\mathrm{T}}S+SBK_s^{k-1}$ to the left side of equation (\ref{lem2pfeq1}) and substitute (\ref{ks}) into it to obtain
\begin{equation*}
\begin{aligned}
&\tilde{A}_{k-1}^{\mathrm{T}}S+S\tilde{A}_{k-1}+\left(K_s^{k-1}\right)^{\mathrm{T}}\Upsilon K_s+K_s^{\mathrm{T}}\Upsilon K_s^{k-1} \\
&-K_s^{\mathrm{T}}\Upsilon K_s-Q_{\Gamma}={\bf0},
\end{aligned}
\end{equation*}
then subtract the above equation from equation (\ref{sk}), giving rise to 
\begin{equation}\label{lem2pfeq3}
\begin{aligned}
\!\!\!&\tilde{A}_{k-1}^{\mathrm{T}}(S_k-S)+(S_k-S )\tilde{A}_{k-1}\\
\!\!\!=&-(K_s-K_s^{k-1})^{\mathrm{T}}\Upsilon(K_s-K_s^{k-1}).
\end{aligned}
\end{equation}
Since $\tilde{A}_{k-1}$ is Hurwitz and $\Upsilon>0$, one has $S_k\geq S$.

In order to prove $\tilde{A}_{k}$ is Hurwitz, equation (\ref{lem2pfeq3}) is transformed into 
\begin{equation*} 
\begin{aligned}
\!\!\!&\tilde{A}_{k}^{\mathrm{T}}(S_k-S)\!+(S_k-S)\tilde{A}_{k}\!\!=\!-(K_s^k\!-\!K_s^{k-1})^{\!\mathrm{T}}\Upsilon\\
\!\!\!&\times (K_s^k-K_s^{k-1})-(K_s-K_s^k)^{\mathrm{T}}\Upsilon(K_s-K_s^k).
\end{aligned}
\end{equation*}
Assume $\tilde{A}_{k}x=\lambda x$ for some $\lambda$ with $\mathrm{Re}(\lambda)\geq0$ and vector $x\neq{\bf0}$. Then, we have
\begin{equation}\label{lem2pfeq5}
\begin{aligned}
&2\lambda x^{\mathrm{T}}(S_k-S)x=-x^{\mathrm{T}}(K_s^{k}-K_s^{k-1})^{\mathrm{T}}\Upsilon(K_s^{k}\\
&-K_s^{k-1})x-x^{\mathrm{T}}(K_s-K_s^k)^{\mathrm{T}}\Upsilon(K_s-K_s^k)x.
\end{aligned}
\end{equation}
As $S_k-S\geq0$, by using 
\begin{equation*}
-\!(\!K_s^{k}\!-\!K_s^{k-1})^{\!\mathrm{T}}\!\Upsilon\!(K_s^{k}\!-\!K_s^{k-1})\!-\!(K_s\!-\!K_s^k)^{\!\mathrm{T}}\!\Upsilon\!(K_s\!-\!K_s^k)\!\leq\!0,
\end{equation*} 
the equality (\ref{lem2pfeq5}) implies 
\begin{equation*}
x^{\mathrm{T}}(K_s^k-K_s^{k-1})^{\mathrm{T}}\Upsilon(K_s^k-K_s^{k-1})x=0,
\end{equation*}
which gives rise to $(K_s^{k}-K_s^{k-1})x={\bf0}$. Consequently, we can get $\tilde{A}_kx=\tilde{A}_{k-1}x=\lambda x$. It contradicts with the induction assumption. 
\end{enumerate}
Therefore, the inequality $S_k\geq S$ and the result 2) hold for all $k\in\mathbb{N}_+$. Next, we show the second inequality of the result 1) also holds for all $k\in\mathbb{N}_+$. 

By using (\ref{ksk}) and completing the square, we rewrite equation (\ref{sk}) as
\begin{equation}\label{sr1}
\begin{aligned}
&\tilde{A}_k^{\mathrm{T}}S_k+S_k\tilde{A}_k-Q_{\Gamma}^k+(K_s^k-K_s^{k-1})^{\mathrm{T}}\\
&\times\Upsilon(K_s^{k}-K_s^{k-1})={\bf0},
\end{aligned}
\end{equation} 
and replace $k$ in (\ref{sk}) by $k+1$, yielding
\begin{equation}\label{sr2}
\tilde{A}_k^{\mathrm{T}}S_{k+1}+S_{k+1}\tilde{A}_k-Q_{\Gamma}^k={\bf0}.
\end{equation}
Subtracting (\ref{sr2}) from (\ref{sr1}), it gives
\begin{equation*}
\begin{aligned}
&\tilde{A}_{k}^{\mathrm{T}}(S_{k}-S_{k+1})+(S_{k}-S_{k+1})\tilde{A}_{k}\\
=&-(K_s^k-K_s^{k-1})^{\mathrm{T}}\Upsilon(K_s^k-K_s^{k-1}),
\end{aligned}
\end{equation*}
by the facts that $\tilde{A}_k$ is Hurwitz and $\Upsilon>0$, which implies that $S_{k+1}\leq S_{k}$, $k\in\mathbb{N}_+$.
 
Therefore, the result 1) is established. It implies that the sequence $\{S_k\}_1^{\infty}$ is monotonically decreasing and has a lower bound $S$. As a result, we conclude that this sequence has a limit. Since $K_s^k$ is the unique solution of (\ref{ksk}), the convergence of $\{K_s^k\}_1^{\infty}$ can be deduced from the convergence of $\{S_k\}_1^{\infty}$. Furthermore, $S$ satisfies equation (\ref{sk}) with $K_s^{k-1}=K_s$, thereby implying that the limit of the sequence $\{S_k\}_1^{\infty}$ is equal to $S$, {\sl i.e., $\lim_{k\rightarrow\infty}S_k=S$}. Consequently, we also have $\lim_{k\rightarrow\infty}K_s^k=K_s$.

Hence, the proof is complete. \hfill $\square$
\end{pf}
Plugging (\ref{kop}) and (\ref{ks}) into (\ref{xbar}), the mean field state dynamics can be written as
% \begin{equation}
% u_i^o=-Kx_i(t)-K_s\bar{x}(t),~~i=1,\cdots,N,
% \end{equation}
% and 
\begin{equation}\label{xbar2}
\dot{\bar{x}}(t)=(A-B(K+K_s))\bar{x}(t),~~\bar{x}(0)=\bar{x}_0,
\end{equation}
where $K$ and $K_s$ are feedback gain matrix and feedforward gain matrix respectively.

In subsequent steps, based on the iterative equations (\ref{kk})-(\ref{pk}) and (\ref{sk})-(\ref{ksk}), as well as the mean field state dynamics (\ref{xbar2}), we aim to further eliminate the need for system dynamics. 
\vspace{-0.5cm}
\subsection{Gain matrix approximation with unknown dynamics}
% \vspace{-1cm}
First, based on iterations of (\ref{pk}) and (\ref{kk}), we develop an equivalent model-free iterative equation to approximate the feedback gain matrix $K$.

To this end, for a fixed $1\leq i\leq N$, introduce the following quadratic function
\begin{equation*}
V_1(x_i)=x_i^{\mathrm{T}}P_kx_i.
\end{equation*}
\vspace{-0.2cm}
By the Ito's formula and SDE (\ref{sys1}), one has
\begin{equation*}
\begin{aligned}
\mathrm{d}V_1= &2x_i^{\mathrm{T}}P_{k}(Ax_i+Bu_i)\mathrm{d}t+(Cx_i+Du_i)^{\mathrm{T}}P_k\\
&\times(Cx_i+Du_i)\mathrm{d}t+2x_i^{\mathrm{T}}P_k(Cx_i+Du_i)\mathrm{d}w_i,
\end{aligned}
\end{equation*}
\vspace{-0.2cm}
which is equivalent to 
\vspace{-0.1cm}
\begin{equation}\label{dv1}
\begin{aligned}
&\!\!\!\!\!\mathrm{d}(x_i^{\!\mathrm{T}}P_kx_i)\!=\!\Big[x_i^{\!\mathrm{T}}(A_{\!k\!-\!1}^{\mathrm{T}}P_k\!+\!P_kA_{\!k\!-\!1}\!+\!C_{\!k\!-\!1}^{\!\mathrm{T}}P_kC_{\!k\!-\!1})x_i\\
&\!\!\!\!\!+\!\!2(u_i\!+\!K_{k-1}x_i)^{\!\mathrm{T}}(B^{\!\mathrm{T}}P_k\!+\!D^{\!\mathrm{T}}P_kC)x_i\!+\!u_i^{\!\mathrm{T}}D^{\!\mathrm{T}}\!P_k\!D\!u_i\!\\
&\!\!\!\!\!-\!x_i^{\!\mathrm{T}}\!K_{k-1}^{\!\mathrm{T}}\!D^{\!\mathrm{T}}\!P_k\!D\!K_{k-1}\!x_i\Big]\!\mathrm{d}t\!+\!2\!x_i^{\!\mathrm{T}}\!P_k\!(Cx_i\!+\!Du_i)\mathrm{d}w_i.
\end{aligned}
\end{equation}
\vspace{-0.2cm}
Substituting (\ref{pk}) and (\ref{kk}) into equation (\ref{dv1}) and letting
\begin{numcases}{}
\Lambda_k=D^{\mathrm{T}}P_kD,\label{ktilde}\\
\tilde{K}_k=(R+\Lambda_k)K_k,\label{lambdak}
\end{numcases}
\vspace{-0.2cm}
yields 
\vspace{-0.1cm}
\begin{equation*} 
\begin{aligned}
&\!\!\!\!\!\!\!\mathrm{d}(x_{i}^{\!\mathrm{T}}P_{k}x_{i})\!\!=\!\!\Big[-x_{\!i}^{\!\mathrm{T}}Q_{\!k\!-\!1\!}x_{i}\!+\!2(u_{i}\!\!+\!\!K_{k\!-\!1}x_{\!i})^{\!\mathrm{T}}\tilde{K}_{\!k}x_{\!i}\!+\!u_{i}^{\!\mathrm{T}}\!\Lambda_k\!u_{i}\\
&\!\!\!\!\!\!-\!\!x_i^{\!\mathrm{T}}\!K_{\!k\!-\!1}^{\mathrm{T}}\Lambda_kK_{ k - 1}x_i \Big]\mathrm{d}t+ 2x_i^{\!\mathrm{T}}P_k(Cx_i\!+\!Du_i)\mathrm{d}w_i.
\end{aligned}
\end{equation*}
Let $t'=t+T$, where $T>0$. We integrate the above equation along SDE (\ref{sys1}) over the time interval $[t,t')$ and obtain
\begin{equation*}
\begin{aligned}
&\!\!x_i(t')^{\!\mathrm{T}}P_{k}x_i(t')\!-\!x_i(t)^{\!\mathrm{T}}P_{\!k}x_i(t)\!=\!-\!\int_t^{t'}\!\!\!x_i^{\!\mathrm{T}}Q_{\!k\!-\!1\!}x_i\mathrm{d}\tau\!\\
&\!\!+\!2\!\int_t^{\!t'}(u_i\!+\!K_{k-1}x_i)^{\!\mathrm{T}}\tilde{K}_{k}x_i\mathrm{d}\tau\!+\!\int_t^{t'}\!\!\!\big(u_i^{\!\mathrm{T}}\Lambda_ku_i-\\
&\!\!x_i^{\mathrm{T}}K_{k\!-\!1}^{\!\mathrm{T}}\Lambda_kK_{k\!-\!1}x_i\big)\mathrm{d}\tau\!+\!2\!\!\int_t^{t'}\!\!\!\!x_i^{\!\mathrm{T}}P_{k}(Cx_i\!\!+\!\!Du_i)\mathrm{d}w_i.
\end{aligned}
\end{equation*}
Then taking the expectation, it gives
\begin{equation}\label{iv1}
\begin{aligned}
\!\!\!&E\left[x_i(t')^{\!\mathrm{T}}P_{k}x_i(t')\!-\!x_i(t)^{\!\mathrm{T}}P_{\!k}x_i(t)\right]\!=\!\!-\!E\bigg[\!\int_t^{t'}x_i^{\!\mathrm{T}}\\
\!\!\!&\times\!Q_{k-1}x_i\mathrm{d}\tau\bigg] \!\!+\!\!2E\bigg[\!\int_t^{\!t'}(u_i\!\!+\!\!K_{k-1}x_i)^{\!\mathrm{T}}\tilde{K}_{k}x_i\mathrm{d}\tau\bigg]\\
&+ E \bigg[ \int_t^{t'}\!\!\!\!\!\big(u_i^{\mathrm{T}}\Lambda_ku_i\!-\!x_i^{\mathrm{T}}K_{k-1}^{\mathrm{T}}\Lambda_kK_{k-1}x_i\big)\mathrm{d}\tau\bigg].\!
\end{aligned}
\end{equation}

Note that the system dynamics have disappeared in equation (\ref{iv1}). To further obtain a data-based equation that can determine the unknown parameters, we define the following $l$-sets-semi-data-based matrices 
\begin{equation}\label{matrx1}
\left\{\begin{aligned}
&\!\Delta_{\hat{x}_i}\!\triangleq\![\delta_{\hat{x}_i}^{t_1},\delta_{\hat{x}_i}^{t_2},\cdots,\delta_{\hat{x}_i}^{t_l}]^{\!\mathrm{T}}, \delta_{\hat{x}_i}^{t}\!\triangleq\!E[\hat{x}_i(t')-\hat{x}_i(t)],\\
&\!\mathcal{I}_{{x}_ix_i}\!\triangleq\![I_{x_ix_i}^{t_1}\!,\!I_{x_ix_i}^{t_2}\!,\!\cdots\!,\!I_{x_ix_i}^{t_l}]^{\!\mathrm{T}}\!,I_{x_ix_i}^{t}\!\!\triangleq\!\!E\!\left[\!\int_t^{t'}\!\!\!\!x_i(\!\tau\!)\!\otimes\! x_i(\!\tau\!)\mathrm{d}\tau\!\right]\!,\\
&\!\mathcal{I}_{x_iu_i}\!\triangleq\![I_{x_iu_i}^{t_1}\!,\!I_{x_iu_i}^{t_2}\!,\!\cdots\!,\!I_{x_iu_i}^{t_l}]^{\!\mathrm{T}}\!,I_{x_iu_i}^{t}\!\!\triangleq\!\!E\!\left[\!\int_t^{t'}\!\!\!\!x_i(\!\tau\!)\!\otimes\! u_i(\!\tau\!)\mathrm{d}\tau\!\right]\!,\\
&\!\mathcal{I}_{\hat{u}_i}\!\triangleq\![I_{\hat{u}_i}^{t_1}\!,\!I_{\hat{u}_i}^{t_2}\!,\!\cdots\!,\!I_{\hat{u}_i}^{t_l}]^{\!\mathrm{T}}\!,I_{\hat{u}_i}^{t}\!\!\triangleq\!\!E\!\left[\!\int_t^{t'}\!\!\!\!\hat{u}_i(\!\tau\!) \mathrm{d}\tau\!\right]\!,
\end{aligned}\right.
\end{equation}
where $l\in\mathbb{N}_{+}$, $t_j=t_1+(j-1)T_s$ represents the time instant, $t_1$ is the initial time, and $T_s>0$ is the sampling time.

Moreover, letting $\rho_j\in\mathbb{R}^n$ be the $j$-th row of $K_{k-1}\in\mathbb{R}^{m\times n}$, we define 
\begin{equation}\label{ktildekminus1}
\!\!\!\!\!\!\left\{\begin{aligned}
&\kappa_{ij}\triangleq \rho_i^{\mathrm{T}}\otimes \rho_j^{\mathrm{T}}\in\mathbb{R}^{n^2}, i,j=1,2,\cdots,m,\\
&\mathbb{K}_{\!k\!-\!1}\!\triangleq\![{\kappa}_{\!11}\!,\!\cdots\!,\!{\kappa}_{\!1m}\!,\!\kappa_{\!22},\!\cdots\!,\!\kappa_{(\!m-1\!)\!m}\!,\!\kappa_{\!mm}]\!\!\in\!\mathbb{R}^{\!n^{\!2}\!\times\!\frac{m}{2}\!(m\!+\!1)}.
\end{aligned}\right.
\end{equation}   

Then, by the Kronecker product and the matrices (\ref{matrx1})-(\ref{ktildekminus1}), we can derive the linear matrix form of equation (\ref{iv1}) as follows
\begin{equation*}
\Psi_{k-1}\left[\begin{array}{c}
\bar{P}_k\\
\mathrm{col}(\tilde{K}_k)\\
\bar{\Lambda}_k
\end{array}\right]=\Xi_{k-1},
\end{equation*}
where  
\begin{equation*} 
\left\{\begin{aligned}
&\Psi_{\!k-\!1}\!=\![\Delta_{\hat{x}_i},\!-\!2\mathcal{I}_{x_iu_i}\!\!-\!2\mathcal{I}_{x_ix_i}\!(\!I\!\otimes\!K_{k-1}^{\mathrm{T}}),\!-\mathcal{I}_{\hat{u}_i}\!+\!\mathcal{I}_{x_ix_i}\!\mathbb{K}_{k-1}],\\
&\Xi_{k\!-1}=-\mathcal{I}_{x_ix_i}\mathrm{col}(Q_{k-1}).
\end{aligned}\right.
\end{equation*}

\begin{assume}\label{A2}
There exists an $l_1>0$ such that for all $l\geq l_1$, we have
\begin{equation}\label{rank1}
\!\!\!\mathrm{rank}\!\!\left(\!\!\left[\begin{array}{cccc}
\!\!I_{x_ix_i}^{\!t_1}\!&\!I_{x_ix_i}^{\!t_2}\!&\!\cdots\!&\!I_{x_ix_i}^{\!t_l}\!\!\\
\!\!I_{x_iu_i}^{\!t_1}\!&\!I_{x_iu_i}^{\!t_2}\!&\!\cdots\!&\!I_{x_iu_i}^{\!t_l}\!\!\\
\!\!I_{\hat{u}_i}^{\!t_1}\!&\!I_{\hat{u}_i}^{\!t_2}\!&\!\cdots\!&\!I_{\hat{u}_i}^{\!t_l}\!\!
\end{array}\right]\!\!\right)\!\!=\!\!\frac{n}{2}
\!(n\!+\!1)\!+\!mn\!+\!\frac{m}{2}\!(m\!+\!1)\!.
\end{equation}
\end{assume}
\begin{remark}\label{rem3.5}
To satisfy the above rank condition mentioned, it is necessary for $l$ that to be not less than $\frac{n}{2}\!(n\!+\!1)\!+\!mn\!+\!\frac{m}{2}\!(m\!+\!1)\!$. There are generally two approaches to fulfill this rank condition. The first approach is collecting data with various control inputs and initial states to ensure that $l$ is sufficiently large. Alternatively, the second approach, which will be used in our algorithm, is introducing an exploration noise into the input channel. This can be achieved by setting $u_i=-Kx_i+\ell_i(t)$, where $K$ is a stabilizer of $[A,B;C,D]$ and $\ell_i(t)$ represents an exploration noise, which is selected from options, such as random noise (\cite{al2007model}), exponentially decreasing probing noise (\cite{vamvoudakis2011multi}), and sum of a series of sinusoidal signals with different frequencies (\cite{jiang2012computational}).
\end{remark}
Now, we are ready to construct a convergent sequence to approximate both the solution of the ARE (\ref{ARE1}) and the feedback gain matrix $K$ as follows.
\begin{theorem}\label{thm1}
Suppose Assumption \ref{A2} holds. Let $P$ be a positive definite solution of (\ref{ARE1}), and let $K$ be a solution of (\ref{kop}). If $K_0\in\mathbb{R}^{m\times n}$ is a stabilizer of $[A,B;C,D]$, then the sequence $\{P_k,\tilde{K}_k,\Lambda_k\}_1^{\infty}$ generated by recursively solving the following equation
\begin{equation}\label{mf1}
\left[\begin{array}{c}
\bar{P}_k\\
\mathrm{col}(\tilde{K}_k)\\
\bar{\Lambda}_k
\end{array}\right]\!=\!\left(\Psi_{k-1}^{\mathrm{T}}\!\Psi_{k-1}\right)^{-1}\!\Psi_{k-1}^{\mathrm{T}}\!\Xi_{k-1},~~k\in\mathbb{N}_+, 
\end{equation}
satisfies the following properties.
\vspace{-0.3cm}
\begin{enumerate}
\item[1).] $\lim_{k\rightarrow\infty}P_k=P$.
\item[2).] $\lim_{k\rightarrow\infty}(R+\Lambda_k)^{-1}\tilde{K}_k=K$.
\end{enumerate} 
\end{theorem}
\vspace{-0.4cm}
\begin{pf}
 Equation (\ref{mf1}) is equivalent to (\ref{pk}) and (\ref{kk}) if the solution of (\ref{mf1}) is unique. Then the convergence result follows Lemma \ref{lem1}. To complete the proof, the remaining task is to show that $\Psi_{k-1}$ has full column rank for all $k\in\mathbb{N}_+$. This will be fulfilled by contradiction.

First, suppose that $\Psi_{k-1} H={\bf0}$ for a nonzero vector $H=[
X_v^{\mathrm{T}},Y_v^{\mathrm{T}},Z_v^{\mathrm{T}}]^{\mathrm{T}}$, where $X_v\in\mathbb{R}^{\frac{n}{2}(n+1)}$, $Y_v\in\mathbb{R}^{mn}$, and $Z_v\in\mathbb{R}^{\frac{m}{2}(m+1)}$. The symmetric matrices $X\in\mathbb{R}^{n\times n}$ and $Z\in\mathbb{R}^{m\times m}$ can be uniquely determined by $\bar{X}=X_v$ and $\bar{Z}=Z_v$, respectively, and a matrix $Y\in\mathbb{R}^{m\times n}$ can be uniquely determined by $\mathrm{col}(Y)=Y_v$. Then, one has
  \begin{equation*}
  \begin{aligned}
&\Delta_{\hat{x}_i}X_v\!-\!\mathcal{I}_{x_ix_i}\mathrm{col}(K_{k-1}^{\mathrm{T}}Y\!+\!Y^{\mathrm{T}}K_{k-1}\!+\!K_{k-1}^{\mathrm{T}}Z\\
&\times K_{k-1})-2\mathcal{I}_{x_iu_i}Y_v-\mathcal{I}_{u_iu_i}\mathrm{col}(Z)={\bf0},
\end{aligned}
\end{equation*}
where 
\begin{equation*}
\mathcal{I}_{u_iu_i}\!\!\!=\!\![\!I_{u_iu_i}^{t_1}\!,\!I_{u_iu_i}^{t_2}\!,\!\cdots\!,\!I_{u_iu_i}^{t_l}\!]^{\!\mathrm{T}},~~I_{u_iu_i}^{t}\!\!\!=\!\!E\!\!\left[\!\int_{t}^{\!t'}\!\!\!u_i\!\otimes \!u_i\mathrm{d}\tau\!\right].
\end{equation*}

On the other hand, by (\ref{dv1}) and (\ref{iv1}), we can get
\begin{equation*}
\begin{aligned}
\!\!\!&\Delta_{\hat{x}_i}X_v \!=\! \mathcal{I}_{x_ix_i}\!\mathrm{col}\Big(\!A_{k-1}^{\!\mathrm{T}}\!X\!+\!X\!A_{k\!-\!1}\!+\!C_{k\!-\!1}^{\!\mathrm{T}}\!X\!C_{k-1}\!+\!\!K_{k-1}^{\!\mathrm{T}}\\
&\!\!\times\!\!(B^{\!\mathrm{T}}\!X\!+\!D^{\!\mathrm{T}}\!X\!C)\!+\!(X\!B\!+\!C^{\!\mathrm{T}}\!X\!D)\!K_{\!k-1}\!\!-\!\!K_{\!k-1}^{\!\mathrm{T}}\!D^{\!\mathrm{T}}\!X\!DK_{\!k\!-\!1}\!\Big)\!\\
&+\!2\mathcal{I}_{x_iu_i}\!\mathrm{col}(B^{\mathrm{T}}X+D^{\mathrm{T}}\!XC_k\!)\!+\!\mathcal{I}_{u_iu_i}\!\mathrm{col}(D^{\mathrm{T}}XD).
\end{aligned}
\end{equation*} 
It follows that
\begin{equation}\label{43}
\begin{aligned}
\mathcal{I}_{\hat{x}_i}\!\bar{M}_1\!+\!2\mathcal{I}_{x_iu_i}\!\bar{M}_2\!+\!\mathcal{I}_{\hat{u}_i}\!\bar{M}_3\!=\!{\bf0},
\end{aligned}
\end{equation}
where $\mathcal{I}_{\hat{x}_i}\!=\![I_{\hat{x}_i}^{t_1}\!,\!I_{\hat{x}_i}^{t_2}\!,\!\cdots\!,\!I_{\hat{x}_i}^{t_l}]^{\!\mathrm{T}}$, $I_{\hat{x}_i}^{t}\!\!=\!\!E\!\left[\!\int_t^{t'}\!\!\!\!\hat{x}_i(\!\tau\!) \mathrm{d}\tau\!\right]$, and 
% $M_1\!\!=\!A_k^{\!\mathrm{T}}X\!+\!X\!A_k\!+\!C_k^{\!\mathrm{T}}\!XC_k\!+\!K_{\!k-1}^{\!\mathrm{T}}M_2+M_2^{\mathrm{T}}K_{k-1}-\!K_{\!k-1}^{\!\mathrm{T}}M_3K_{\!k-1}$, $M_2=B^{\mathrm{T}}X+D^{\mathrm{T}}XC_k-Y$, and $M_3=D^{\mathrm{T}}XD-Z$.
\begin{equation*}
\left\{\begin{aligned}
&M_1=A_{k-1}^{\mathrm{T}}X\!+\!X\!A_{k-1}\!+\!C_{k-1}^{\mathrm{T}}\!X\!C_{k-1}\!+\!K_{k-1}^{\mathrm{T}}M_2\!\\
&~~~~~~+\!M_2^{\mathrm{T}}K_{k-1}\!-\!K_{\!k-1}^{\!\mathrm{T}}M_3K_{\!k-1},\\
&M_2=B^{\mathrm{T}}X+D^{\mathrm{T}}XC-Y,\\
&M_3=D^{\mathrm{T}}XD-Z.
\end{aligned}\right.
\end{equation*}
Under the rank condition (\ref{rank1}) of Assumption \ref{A2}, equation (\ref{43}) gives rise to $M_1={\bf0}$, $M_2={\bf0}$, and $M_3={\bf0}$. As a result, we can get
\begin{equation}
A_{k-1}^{\mathrm{T}}X+XA_{k-1}+C_{k-1}^{\mathrm{T}}XC_{k-1}={\bf0},
\end{equation}
followed with the fact that $K_{k-1}$ is a stabilizer of $[A,B;C,D]$ and \cite[Theorem 1]{rami2000linear}, which implies $X={\bf0}$. It further gives $Y={\bf0}$ and $Z={\bf0}$, sequentially.

Thus, we have $H={\bf0}$, leading to a contradiction with the assumption that $X\neq{\bf0}$, which states that $\Psi_{k-1}$ must have full column rank for all $k\in\mathbb{N}_+$. 

Hence, the proof is complete. \hfill $\square$
\end{pf}
\begin{remark}
In reference \cite{li2022stochastic}, the proposed method successfully removes the need for the state matrix $A$. However, it still relies on  the knowledge of the input-to-state matrix $B$ and diffusion coefficients $C,D$. In contrast, our method goes to a step further by completely eliminating the system dynamics from the iterative equation. Moreover, the approximation of $D^{\mathrm{T}}PD$ obtained by the above algorithm can be further used in developing a model-free iterative equation to approximate $(S,K_s)$, as we introduced later.
\end{remark}
Second, we continue to eliminate the system information in the iterative equations (\ref{sk})-(\ref{ksk}). Rewrite the dynamics of the $i$-th agent as
\begin{equation}
\mathrm{d}x_i\!=\!\Big[\!\tilde{A}_{\!k\!-\!1}\!{x}_i\!+\!B\!\big(\!(\!K\!+\!K_s^{\!k-1}\!)x_i\!+\!u_i\!\big)\!\Big]\!\mathrm{d}t\!+\!(\!Cx_i\!+\!Du_i\!)\mathrm{d}w_i,
\end{equation}
and take the expectation on both sides of the above equation to obtain
\begin{equation}\label{sysxibar0}
\mathrm{d}\bar{x}_i=\Big[\tilde{A}_{k-1}\bar{x}_i\!+\!B\big((K+K_s^{k-1})\bar{x}_i+\bar{u}_i\big)\Big]\mathrm{d}t,
\end{equation}
where $\bar{x}_i=E[x_i]$ and $\bar{u}_i=E[u_i]$.

Then, taking the integration operation for $\frac{\mathrm{d}}{\mathrm{d}t}(\bar{x}_iS_k\bar{x}_i)$ and using equations  (\ref{sysxibar0}) and (\ref{sk})-(\ref{ksk}), we can deduce that
\begin{equation}\label{sksks}
\begin{aligned}
\!\!\!&\bar{x}_i(t')^{\!\mathrm{T}}S_k\bar{x}_i(t')\!-\!\bar{x}_i(t)^{\!\mathrm{T}}S_k\bar{x}_i(t)\!-\!2\!\!\int_{t}^{t'}\!\!\!\!\big(\!(K_s^{k-1}\!+\!K)\bar{x}_i\\
\!\!\!&~~+\!\bar{u}_i\big)^{\!\mathrm{T}}\Upsilon{K}_s^k\bar{x}_i\mathrm{d}\tau+\!\int_{t}^{t'}\!\bar{x}_i^{\mathrm{T}}{Q}_{\Gamma}^{k-1}\bar{x}_i\mathrm{d}\tau\!=\!0.
\end{aligned}
\end{equation}
To further develop an iterative equation based on this model-free equation, we define new $l$-sets-data-based matrices 
\begin{equation}\label{matrix2}
\left\{\begin{aligned}
&\Delta_{\hat{\bar{x}}_i}\triangleq[\delta_{\hat{\bar{x}}_i}^{t_1},\delta_{\hat{\bar{x}}_i}^{t_2},\cdots,\delta_{\hat{\bar{x}}_i}^{t_l}]^{\!\mathrm{T}},~~\delta_{\hat{\bar{x}}}^{t}\triangleq\hat{\bar{x}}_i(t')-\hat{\bar{x}}_i(t),\\
&\mathcal{I}_{\bar{x}_i\bar{x}_i}\!\triangleq\![I_{\bar{x}_i\bar{x}_i}^{t_1},\!I_{\bar{x}_i\bar{x}_i}^{t_2},\!\cdots,\!I_{\bar{x}_i\bar{x}_i}^{t_l}]^{\!\mathrm{T}}\!,I_{\bar{x}_i\bar{x}_i}^{t}\!\triangleq\!\!\int_t^{t'}\!\!\!\!\bar{x}_i(\!\tau\!)\!\otimes\! \bar{x}_i(\!\tau\!)\mathrm{d}\tau,\\
&\mathcal{I}_{\bar{x}_i\bar{u}_i}\!\triangleq\![I_{\bar{x}_i\bar{u}_i}^{t_1},\!I_{\bar{x}_i\bar{u}_i}^{t_2},\!\cdots,\!I_{\bar{x}_i\bar{u}_i}^{t_l}]^{\!\mathrm{T}}\!,I_{\bar{x}_i\bar{u}_i}^{t}\!\triangleq\!\!\int_t^{t'}\!\!\!\!\bar{x}_i(\!\tau\!)\!\otimes\! \bar{u}_i(\!\tau\!)\mathrm{d}\tau,
\end{aligned}\right.
\end{equation}
and then derive the linear matrix of (\ref{sksks}) as follows
\begin{equation}
\Phi_{k-1}\left[\begin{array}{c}
\bar{S}_k\\
\mathrm{col}({K}_s^k)
\end{array}\right]= \Theta_{k-1},
\end{equation} 
where 
\begin{equation}\label{par2}
\left\{\begin{aligned}
&\Phi_{k-1}\!\!=\!\!\left[\Delta_{\hat{\bar{x}}_i}\!,\!-\!2\mathcal{I}_{\bar{x}_i\bar{u}_i}\!(\!I\!\otimes\!\Upsilon)\!-\!2\mathcal{I}_{\bar{x}_i\bar{x}_i}\!\!\left(\!I\!\otimes\!\left(\Upsilon\!(K_s^{\!k-1}\!\!+\!K)\right)^{\!\!\mathrm{T}}\!\right)\!\right]\!,\\
&\Theta_{k-1}=-\mathcal{I}_{\bar{x}_i\bar{x}_i}\mathrm{col}({Q}_{\Gamma}^{k-1}).
\end{aligned}\right.
\end{equation}

\begin{assume}\label{A3}
There exists an $l_2>0$ such that for all $l\geq l_2$, we have
\begin{equation}\label{rank2}
\!\!\!\mathrm{rank}\!\!\left(\left[\begin{array}{cccc}
\!\!I_{\bar{x}_i\bar{x}_i}^{\!t_1}\!&\!I_{\bar{x}_i\bar{x}_i}^{\!t_2}\!&\!\cdots\!&\!I_{\bar{x}_i\bar{x}_i}^{\!t_l}\!\!\\
\!\!I_{\bar{x}_i\bar{u}_i}^{\!t_1}\!&\!I_{\bar{x}_i\bar{u}_i}^{\!t_2}\!&\!\cdots\!&\!I_{\bar{x}_i\bar{u}_i}^{\!t_l}\!\!
\end{array}\right] \right)\!\!=\!\!\frac{n}{2}
\!(n\!+\!1)\!+\!mn.
\end{equation}
\end{assume}
\begin{remark}
To meet the above rank condition, it is necessary for $l$ to exceed $\frac{n}{2}(n+1)+mn$. This requirement can be satisfied by the introduction of diverse exploration noise, as mentioned in Remark \ref{rem3.5}. Under the condition, we can design a model-free algorithm that can approximate both the solution of the ARE (\ref{ARE2-2}) and the feedforward gain matrix $K_s$ as follows.
\end{remark}
\begin{theorem}\label{thm2}
Suppose Assumption \ref{A3} holds. Let $S$ be a maximal solution of (\ref{ARE2-2}), and let $K_s$ be a solution of (\ref{ks}). If $K_s^0={\bf0}$, then the sequence $\{S_k,K_s^k\}_1^{\infty}$ generated by recursively solving the following equation
\begin{equation}\label{mf2}
\left[\begin{array}{c}
\bar{S}_k\\
\mathrm{col}({K}_s^k)
\end{array}\right]= \left(\Phi_{k-1}^{\mathrm{T}}\Phi_{k-1}\right)^{\mathrm{T}}\Phi_{k-1}^{-1}\Theta_{k-1},~k\in\mathbb{N}_+,
\end{equation}
with $\Phi_{k-1}$ and $\Theta_{k-1}$ given in (\ref{par2})  satisfies the following properties.
\vspace{-0.3cm}
\begin{enumerate}
\item[1).] $\lim_{k\rightarrow\infty}S_k=S$.
\item[2).] $\lim_{k\rightarrow\infty}K_s^k=K_s$. 
\end{enumerate}
\end{theorem}
\bproof
If the solution of equation (\ref{mf2}) is unique, it can be established that equation (\ref{mf2}) is equivalent to equations (\ref{sk}) and (\ref{ksk}). Remark \ref{rem3.3} states that $A-BK$ is Hurwitz, which implies that the matrix $A-B(K+K_s^0)$ has eigenvalues with negative real parts. As a result, according to \cite[Lemma 6]{jiang2012computational}, the full column rank condition of $\Psi_{k-1}$ for all $k\in\mathbb{N}_+$ is ensured, consequently guaranteeing the uniqueness of the solution of equation (\ref{mf2}). Therefore, the convergence result of $\{S_k,K_s^k\}$ follows with Lemma \ref{lem3.2}.
% Since $A-BK$ is Hurwitz as mentioned in Remark \ref{rem3.3}, which gives that $A-B(K+K_s^0)$ is Hurwitz. As a results, within the intial condition $K_s^0={\bf0}$, the full column rank condition of $\Psi_k$ for all $k\in\mathbb{N}_+$ guranteed by \cite[Lemma 6]{jiang2012computational} ensure the uniqueness of the solution to (\ref{mf2}). Then, the convergence result of the sequence $\{S_k,K_s^k\}_1^{\infty}$ follows with Lemma \ref{lem3.2}.
\eproof

\subsection{Mean field state computation with unknown dynamics}
Finally, to compute the mean field state $\bar{x}(t)$, we substitute
\begin{equation}\label{uxbar}
u_i=-(K+K_s)x_i 
\end{equation}
into equation (\ref{sys1}), yielding
\begin{equation}\label{sysxibar}
\mathrm{d}x_i\!=\!(A\!-\!B(K\!+\!K_s))x_i\mathrm{d}t\!+\!(C-D(K+K_s))x_i\mathrm{d}w_i,
\end{equation}
and set the initial state $x_i(0)=\bar{x}_0$. Then, combining (\ref{xbar2}) and (\ref{sysxibar}), we can get
\begin{equation*} 
\begin{aligned}
\bar{x}(t)=E[x_i(t)].
\end{aligned}
\end{equation*}
\begin{remark}
In practical implementations, the mean field state can be approximated by taking the average of state samples collected from a fixed agent $\mathcal{A}_i$, which is driven by (\ref{uxbar}). Let $x_i^j$ represent the $j$-th state sample of $\mathcal{A}_i$ and $N_s$ be the number of samples. By the law of large numbers, the average of these state samples converges to the MF state as $N_s\rightarrow\infty$ {\sl i.e.,$\lim_{N_s\rightarrow\infty}\frac{1}{N_s}\sum_{j=1}^{N_s}x_i^j(t)=\bar{x}(t)$.}
\end{remark}

In addition, when $S$ is invertible, $S_k$ is also invertible for sufficiently large values of $k$. In this case, the mean field state can be computed using the collected data (\ref{matrix2}) via identifying the system matrices.

By using the final values of the sequence $\{S_k,K_s^k\}$, $\Upsilon$, and equation (\ref{ksk}), the system matrix $B$ can be calculated by
\begin{equation}\label{Bc}
B=(\Upsilon K_s^kS_k^{-1})^{\mathrm{T}}.
\end{equation}
To identify the matrix $A$, we rewrite (\ref{sysxibar0}) as 
\begin{equation}\label{sysxibar2}
\frac{\mathrm{d}\bar{x}_i}{\mathrm{d}t}= A\bar{x}_i+B\bar{u}_i,
\end{equation}
In addition, let $[A]_j$ be the $j$-th row of the matrix $A$, and introduce the quadratic function $\bar{x}_i(t)^{\mathrm{T}}E_j\bar{x}_i(t)$ and differentiate it with respect to time to obtain
\begin{equation*}
\frac{\mathrm{d}}{\mathrm{d}t}(\bar{x}_i^{\mathrm{T}}E_j\bar{x}_i)=2\bar{x}_i^{\mathrm{T}}e_j^{\mathrm{T}}\left([A]_j\bar{x}_i+e_jB\bar{u}_i\right).
\end{equation*}
By integrating the above equation along (\ref{sysxibar2}) over the time interval $[t,t')$ and using the Kronecker product representation, we obtain 
\begin{equation}\label{aaaa}
\!\!\!(\!\delta_{\hat{\bar{x}}_i}^t\!)^{\!\mathrm{T}}\!\bar{E}_j\!=\!2(\!I_{\bar{x}_i\!\bar{x}_i}^t\!)^{\!\mathrm{T}}\!(\!e_j^{\!\mathrm{T}}\!\otimes\!I) [A]_{j}^{\mathrm{T}}\!\!+\!2(\!I_{\!\bar{x}_i\!\bar{u}_i}^{t}\!)^{\!\mathrm{T}}\!\mathrm{col}(\!B^{\!\mathrm{T}}\!E_j\!).
\end{equation}
Then, using the $l$-sets-data-based matrices (\ref{matrix2}), equation (\ref{aaaa}) can be expressed by 
\begin{equation} 
Z_j[A]_j^{\mathrm{T}}\!=H_j,
\end{equation}
where
\begin{equation*}
\!\!\!\!\left\{\begin{aligned}
&Z_j=2\mathcal{I}_{\bar{x}_i\bar{x}_i}\!( e_j^{ \mathrm{T}}\otimes \!I ),\\
&H_j= \Delta_{\hat{\bar{x}}_i}\!\bar{E}_j-2\mathcal{I}_{\bar{x}_i\bar{u}_i}\!\mathrm{col}(B^{\!\mathrm{T}}\!E_j),~~j=1,2,\cdots,n.
\end{aligned}\right.
\end{equation*}
Under Assumption \ref{A3}, it is easy to show that the rank condition $\mathrm{rank}(\mathcal{I}_{\bar{x}_i\bar{x}_i}\!(\!e_j^{\!\mathrm{T}}\!\otimes I ))=n$ always holds for all $j=1,2,\cdots,n$, which gives that each row of the matrix $A$ can be calculated by
\begin{equation}\label{Ac}
\begin{aligned}
[A]_j^{\mathrm{T}}\!=\!\left(Z_j^{\mathrm{T}}Z_j\right)^{-1}\!\left(Z_j^{\mathrm{T}}H_j\right),~~j=1,2,\cdots,n.
\end{aligned}
\end{equation}
Thus, the mean field state can be calculated as
\begin{equation}\label{xbarc}
\bar{x}(t)=e^{(A-BK-BK_s)t}\bar{x}_0,~~t\geq0.
\end{equation}

Herein, we are in the position to present the  model-free mean field social optimal control design algorithm.

\begin{algorithm}[!htb]
\caption{Model-free mean field social optimal control design}\label{alg1}
\begin{algorithmic}
\State {\bf{[Initialization]}}: Choose a stabilizer $K_0$ of the system $[A,B;C,D]$ and $K_s^0={\bf0}$. Select $\mathcal{A}_1$ for learning. Set $k=1$ and convergence criterion $\xi$. 
\State {\bf{[Data Collection]}}: Employ $u_1=-K_0x_1+\ell_1$ to $\mathcal{A}_1$ on $[t_1,t_l]$, where $\ell_1$ is the exploration noise, and other agents remain in place. Along state and input trajectory samples from $\mathcal{A}_1$, calculate matrices (\ref{matrx1}) and (\ref{matrix2}) until conditions (\ref{rank1}) and (\ref{rank2}) are satisfied.

\State [{\bf Approximation of $(P,K,D^{\mathrm{T}}PD)$}]:  Repeatedly solve $(P_k,\tilde{K}_k,\Lambda_k)$ from equation (\ref{mf1}) and let $k=k+1$ until $\|K_{k}-K_{k-1}\|\!\leq\!\xi$. Then, let $\hat{P}=P_k$, $\widehat{D^{\mathrm{T}}PD}\!=\!\Lambda_k$, and $\hat{K}\!=\!(\Lambda_k\!+\!R)^{\!-1}\tilde{K}_k$. Reset $k=1$.

\State [{\bf Approximation of $(S,K_s)$}]: Replace $\Upsilon$ with $(R+\widehat{D^{\mathrm{T}}PD})$ in (\ref{mf2}). Repeatedly solve $(S_k,K_s^k)$ from equation (\ref{mf2}) and let $k=k+1$ until $\|K_s^{k}-K_s^{k-1}\|\leq \xi$, then let $\hat{S}=S^k$ and $\hat{K}_s=K_s^{k}$.
\State [{\bf Approximation of $\bar{x}(t)$ \uppercase\expandafter{\romannumeral1} }]: Employ $u_1=-(\hat{K}+\hat{K}_s)x_1$ to the selected $\mathcal{A}_1$, set $x_1(0)\!=\!\bar{x}_0$, and collect state trajectory samples $x_1^j(t)$, $j\!=\!1\!,\!2\!,\!\cdots,N_s$. Compute $\hat{\bar{x}}(t)\!=\!\frac{1}{N_s}\sum_{j=1}^{N_s}x_1^j$.
\State [{\bf Approximation of $\bar{x}(t)$ \uppercase\expandafter{\romannumeral2} }]: Identify the matrices $B$ and $A$ from (\ref{Bc}) and (\ref{Ac}) and then compute the approximate mean field state $\hat{\bar{x}}(t)$ from (\ref{xbarc}).
\State [{\bf Output}]: Apply $u_i=-\hat{K}x_i(t)-\hat{K}_s\hat{\bar{x}}(t)$ to $N$ agents.
\end{algorithmic}
\end{algorithm}

The effectiveness of the algorithm can be supported by theoretical guarantees provided by Theorems \ref{thm1} and \ref{thm2}.
\begin{figure*}[!htb]
\centering
\includegraphics[scale=0.78]{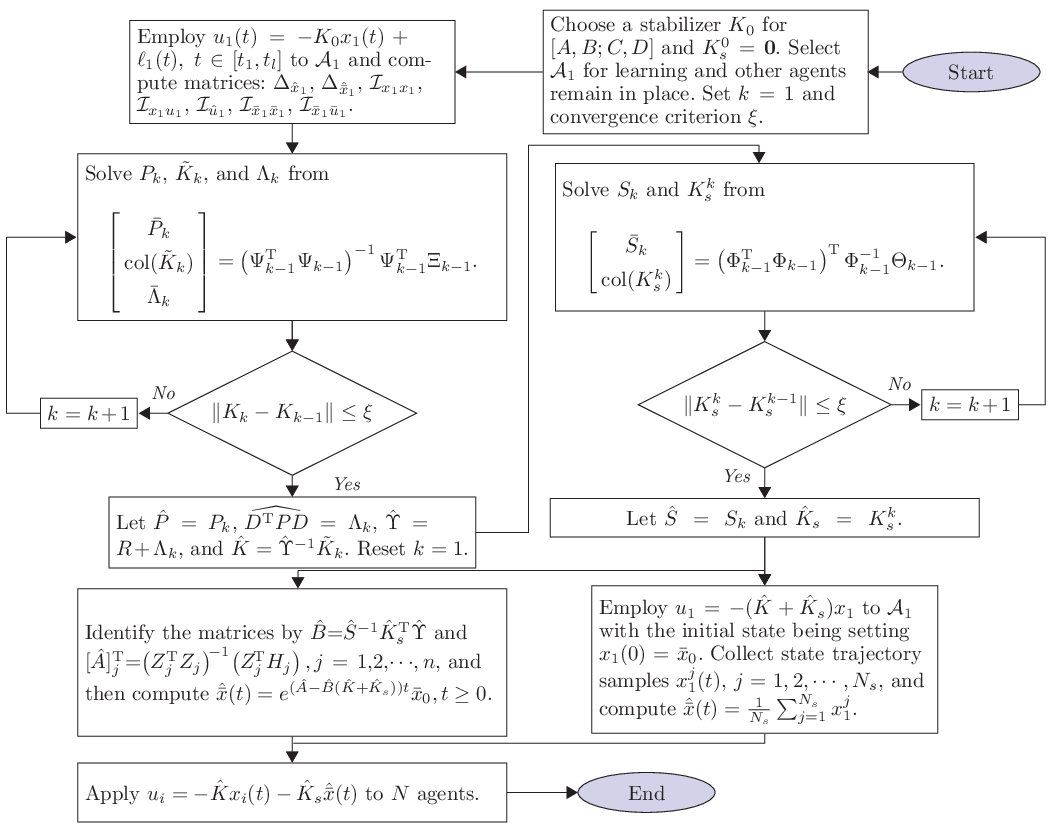}
\centering
\caption{Flowchart of Algorithm \ref{alg1}.}\label{flowchart}
\end{figure*}
\begin{remark}
The flowchart of Algorithm \ref{alg1} is given in Fig. \ref{flowchart}. In this algorithm, the gain matrices associated with the optimal MF social control policy are updated in a sequential manner, followed by a one-time calculation of the MF state. Notably, the two policy iterations characterize the off-policy nature, allowing us to utilize previously generated data to learn the gain matrices online. As a result, the proposed model-free method eliminates the need for acquiring new samples in data-driven iterations, which effectively reduces the computational burden and significantly enhances the efficiency of the algorithm. 
% Denote the end of learning time is $T_l$, the 
\end{remark}
\begin{remark}
Note that during the time interval $[t_1,t_l]$, the control input is exclusively applied to a specified agent, enabling the collection of state and input samples. Once the rank conditions (\ref{rank1}) and (\ref{rank2}) are met, the off-policy RL iterations of $\{P_k,\tilde{K}_k,\Lambda_k\}$ and $\{S_k,K_s^k\}$ are performed sequentially. After this, the mean field state $\bar{x}(t)$ is computed offline. Finally, the obtained control policy is applied to all agents.
\end{remark}

\section{Simulation Example}\label{sec4}
In this section, a numerical simulation is carried out to validate the effectiveness of the proposed algorithm. The large-scale population involves $40$ agents, where coefficients of each agent's dynamics are as follows
\begin{equation*} 
\begin{aligned}
A\!=\!\left[\!\begin{array}{cc}
0.3&0.7\\
-0.9&0.5
\end{array}\!\right]\!,\!B\!=\!\left[\!\begin{array}{c}
0.2\\
0
\end{array}\!\right]\!,\!C\!=\!\left[\!\begin{array}{cc}
0.05&0.03\\
0.05&0.02
\end{array}\!\right]\!,\!D\!=\!\left[\!\begin{array}{cc}
0.05\\
0.06
\end{array}\!\right]\!,
\end{aligned}
\end{equation*}
$x_i\in\mathbb{R}^2$, $u\in\mathbb{R}$, and $w_i$ is a standard one-dimensional Brownian motion. The initial state $x_i(0)$ is uniformly distributed on $[0,4]\times[0,4]\subset\mathbb{R}^2$ with $Ex_i(0)=[2,2]^{\mathrm{T}}$. 

The coefficients of the cost function (\ref{Jsoc}) are 
\begin{equation*}
\begin{aligned}
Q=\left[\begin{array}{cc}
3&0\\
0&2
\end{array}\right],~\Gamma=\left[\begin{array}{cc}
0.9&0\\
0&0.9
\end{array}\right],~R=1.25.
\end{aligned}
\end{equation*}

% which give rise to the following gain matrices 
% \begin{equation*}
% \begin{aligned}
% \!\!\!K_1^*\!=\!\left[\!\begin{array}{cc}
% \!56.9890 & 34.0424
% \end{array}\!\right],K_2^*\!=\!\left[\!\begin{array}{cc}
% \!-20.1499 & 2.7483
% \end{array}\!\right].
% \end{aligned}
% \end{equation*} 
% Clearly, Assumptions ({\bf A1})-({\bf A3}) are verifiable.

In this simulation, to implement Algorithm \ref{alg1}, the first agent $\mathcal{A}_1$ is selected, and the control input of $\mathcal{A}_1$ and is designed as
\begin{equation}\label{u12}
\begin{aligned}
K_0\!=\![6,-3],~~\ell_1(t)\!=\!\sum_{j=1}^{100}\!\sin(w_jt),
\end{aligned}
\end{equation}
where $w_j$, $j=1,\cdots,100$ are randomly selected from $[-100,100]$. 

The other parameters are set to $T_s=0.001$ [sec], $T=0.9$ [sec], the learning duration $[t_1,t_l]=[0,10]\text{~[sec]}$, and the convergence criterion $\xi=10^{-4}$.

\begin{figure*}[!htb]
\centering
\includegraphics[scale=0.33]{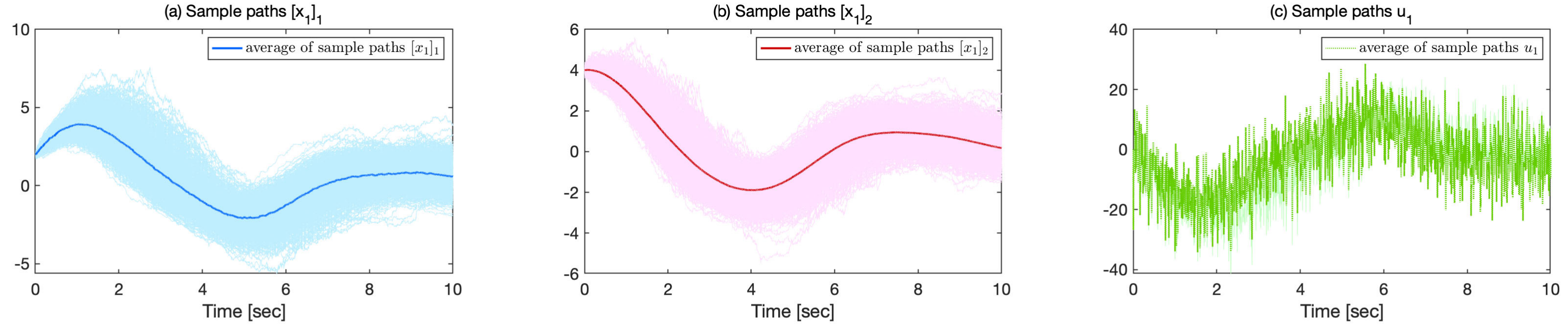}
\centering
\caption{Real-time data collected from $\mathcal{A}_1$.}\label{rf1}
\end{figure*}

To achieve the rank conditions (\ref{rank1}) and (\ref{rank2}), we employ the control input  (\ref{u12}) to drive $\mathcal{A}_1$. Subsequently, data is collected from $100$ sample paths to compute the data-based matrices. Fig. \ref{rf1} shows all $100$ state sample paths and input sample paths for $\mathcal{A}_1$ driven by (\ref{u12}), along with their respective averages.

\begin{figure*}[!htb]
\centering
\includegraphics[scale=0.33]{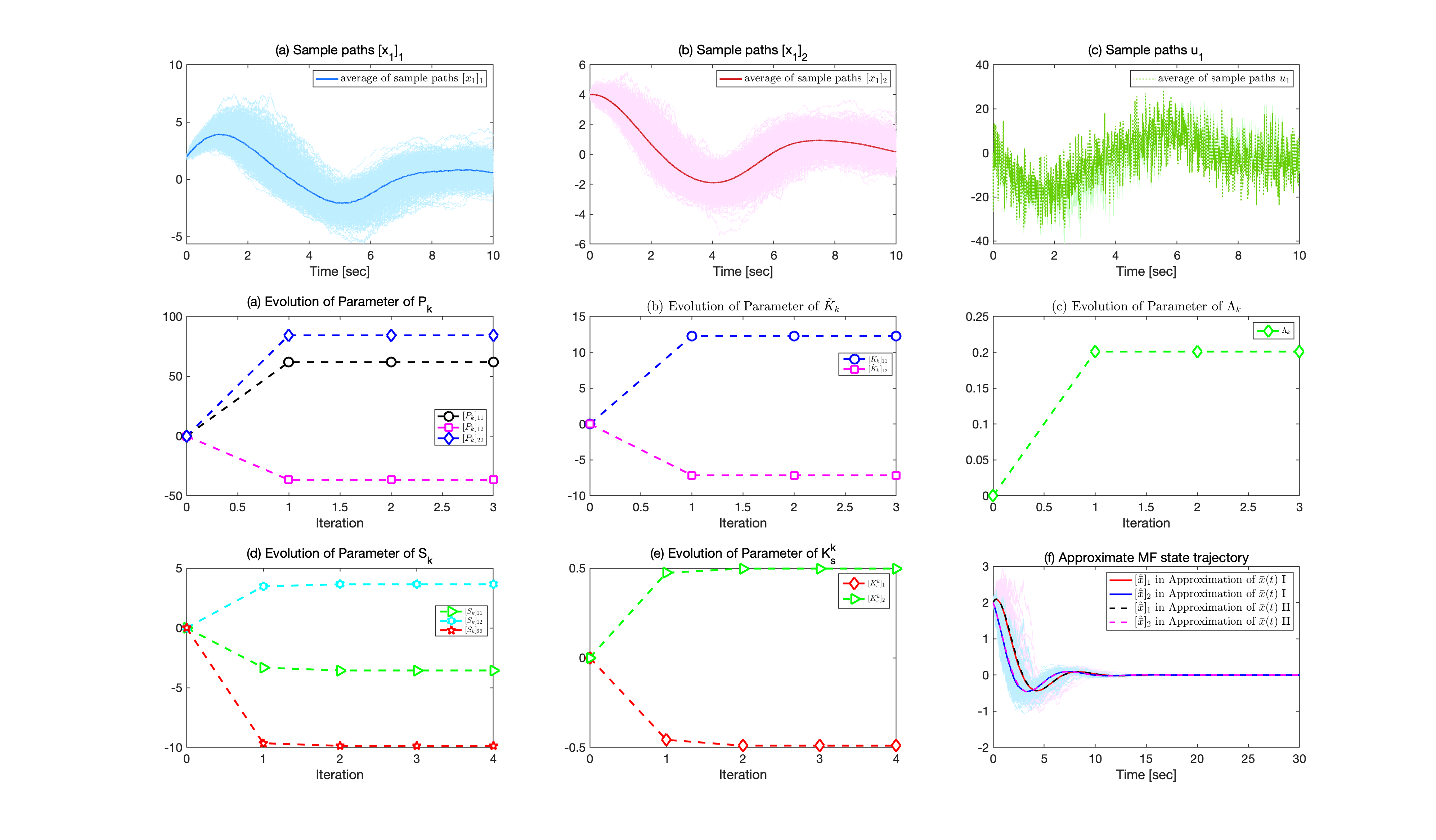}
\centering
\caption{Determined parameters of Algorithm \ref{alg1}.}\label{rf2}
\end{figure*}

The results of equation (\ref{mf1}) in Algorithm \ref{alg1} are presented in Fig. \ref{rf2} (a)-(c), where we observe that the sequence $\{P_k,\tilde{K}_k,\Lambda_k\}$ converges at the $3$rd iteration , as per our chosen convergence criterion. Tab. \ref{e1tab1} summarizes the final values of the estimated parameters and their corresponding approximate errors, where 
\begin{equation*}
\left\{\begin{aligned}
&\mathcal{R}(\hat{P})=A^{\mathrm{T}}\hat{P}+\hat{P}A+C^{\mathrm{T}}\hat{P}C-(B^{\mathrm{T}}\hat{P}+D^{\mathrm{T}}\hat{P}C)^{\mathrm{T}}\Upsilon^{-1}\\
&~~~~~\times(B^{\mathrm{T}}\hat{P}+D^{\mathrm{T}}\hat{P}C)+Q,\\
&\mathcal{K}(\hat{P})=(R+D^{\mathrm{T}}\hat{P}D)^{-1}(B^{\mathrm{T}}\hat{P}+D^{\mathrm{T}}\hat{P}C),\\
&\mathcal{D}(\hat{P})=D^{\mathrm{T}}\hat{P}D.
\end{aligned}\right.
\end{equation*}
It indicates that the small value of $\|\mathcal{R}(\hat{P})\|_2$ confirms the accuracy of $\hat{P}$ as an estimate for the real solution of (\ref{ARE1}). Furthermore, the negligible values of the other errors suggest that both $\hat{K}$ and $\widehat{D^{\mathrm{T}}PD}$ are close to their actual values.
\begin{table}[!htb]
\begin{center}
\begin{tabular}{|c|c|c|c|}
\hline
Parameter & Value & \multicolumn{1}{c|}{Error}             & Value             \\ \hline
$[\hat{P}]_{11}$ & $61.8000$   & \multirow{3}{*}{$\|\mathcal{R}(\hat{P})\|_2$} & \multirow{3}{*}{$0.7299$} \\ 
$[\hat{P}]_{12}$ & $-36.5983$  &                   &                   \\ 
$[\hat{P}]_{22}$ & $84.2412$   &                   &                   \\ \hline
$[\hat{K}]_{11}$ & $8.4670$    & \multirow{2}{*}{$\|\hat{K}-\mathcal{K}(\hat{P})\|_2$} & \multirow{2}{*}{$0.0698$} \\ 
$[\hat{K}]_{12}$ & $-4.9231$   &                   &                   \\ \hline
$\widehat{D^{\mathrm{T}}PD}$ & $0.2010$ &  $|\widehat{D^{\mathrm{T}}PD}-\mathcal{D}(\hat{P})|$              &       $ 0.0371$            \\ \hline
\end{tabular}
\end{center}
\caption{Estimates of $P$, $K$, and $D^{\mathrm{T}}PD$.}\label{e1tab1}
\end{table}

Based on the above results, the sequence $\{S_k,K_s^k\}$ can be obtained by iteratively solving (\ref{mf2}), and the convergence results are shown in Fig. \ref{rf2}. (d) and (e). The simulation demonstrates that $\{P_{11}^k,K_1^k\}$ converges at the $4$th iteration, and Tab. \ref{e1tab2} lists the final values of the determined parameters and their true values. These results imply that $S_k$ and $K_s^k$ of Algorithm \ref{alg1} can converge to their true values.
\begin{table}[!htb]
\begin{center}
\begin{tabular}{|c|c|c|c|}
\hline
Parameter & Value & Parameter & Value \\ \hline
$[\hat{S}]_{11}$ &  $-3.5591$ & $[{S}]_{11}$  &   $-3.4935$    \\ 
$[\hat{S}]_{12}$ &  $3.6498$  & $[{S}]_{12}$ &    $3.5718$   \\ 
$[\hat{S}]_{22}$ &  $-9.8665$ & $[{S}]_{22}$ & $-9.7025$      \\ \hline
$[\hat{K}]_{11}$ &  $-0.4899$ & $[{K}]_{11}$  &   $-0.4815$    \\ 
$[\hat{K}]_{12}$ &  $0.4977$ &  $[{K}]_{12}$  &   $0.4923$    \\ \hline
\end{tabular}
\end{center}
\caption{Estimates of $S$ and $K_s$.}\label{e1tab2}
\end{table}

To approximate the mean field state trajectory, we reset the initial condition of $\mathcal{A}_1$ to $x_1(0)=[2,2]^{\mathrm{T}}$ and enforce $u_1=-(\hat{K}+\hat{K}_s)x_1$ to $\mathcal{A}_1$. We obtain the mean field state trajectory by averaging the state sample paths of $\mathcal{A}_1$. The result is shown illustrated in Fig. \ref{rf2} (f).

On the other hand, by using the previously collected data, we calculate the system matrices by using (\ref{Bc}) and (\ref{Ac}), the identified results are as follows
\begin{equation*}
\hat{A}=\left[\begin{array}{cc}
0.3028&0.7082\\
-0.8887&0.4995
\end{array}\right],~~\hat{B}=\left[\begin{array}{c}
0.2009\\
0.0011
\end{array}\right].
\end{equation*}
and then the the mean field state is approximated by
\begin{equation*}
\hat{\bar{x}}(t)= \mathrm{exp}\left(\left[\begin{array}{cc}
-1.2995& 1.5971\\
   -0.8975 & 0.5044
\end{array}\right]t\right)\left[\begin{array}{c}
2\\
2
\end{array}\right].
\end{equation*}
This estimated mean field state trajectory is also shown in Fig. \ref{rf2} (f), demonstrating the agreement of the mean field state approximations obtained from the aforementioned two methods.

\section{Conclusion}\label{sec5}
In the context of solving MFG problems, most model-free methods rely on a double-loop structure, where an inner loop determines the optimal control policy given the current mean field state, and an outer loop updates the mean-field state. However, as shown in this paper, in the LQG setting, simpler optimality conditions can often be expressed in terms of the Riccati equations governing gain matrices, as well as ordinary or stochastic differential equations describing the evolution of the mean field state. This fact enabled us to obtain the agent’ s own state feedback gain by data-driven policy iteration, which leads to further derive the mean field feedforward gain by technically introducing an equation transformation and using off-policy algorithm. Finally, the mean field state can be approximated off-line. It should be emphasized that the proposed policy iterations can also serve as computational tool for LQG stochastic optimal control problems.

          % and a bib file to produce the 

% \section{Some Latin vocabulary}         % Sections and subsections are supported  
\appendix                         % in the appendices.
\section{Proof of Lemma \ref{lem1}}\label{pf0}

The convergence of $\{K_k\}_1^{\infty}$ is guaranteed by the convergence of $\{P_{k}\}_{1}^{\infty}$ as indicated by (\ref{kk}). To complete the proof, we must establish $\lim_{k\rightarrow\infty}P_k=P$. This is achieved by demonstrating the following results hold for all $k\in\mathbb{N}_+$.
\begin{enumerate}
\item[1).] $P_k$ is a unique positive definite solution of (\ref{pk}).
\item[2).] $K_{k}$ is a stabilizer of $[A,B;C,D]$.
\item[3).] $P_k\geq P_{k+1}\geq P$.
\end{enumerate}

Suppose that $K_{k-1}$ is a stabilizer of $[A,B;C,D]$ for $k\geq1$, we first prove by contradiction that $[A_{k-1},\!C_{k-1};\!\sqrt{Q_{k-1}}]$ is exactly observable. Suppose otherwise; then, by the stochastic Popov-Belevith-Hautus criterion for exact observability (\cite[Theorem 4]{zhang2004stabilizability}), there exists a nonzero symmetric matrix $X$ such that
  \begin{equation}\label{pf1eq1}
  \left\{\begin{aligned}
  &XA_{k-1}^{\mathrm{T}}\!+\!A_{k-1}X\!+\!C_{k-1}XC_{k-1}^{\mathrm{T}} \!=\! \lambda X, \\
  &\sqrt{Q_{k-1}}X = {\bf 0 },~~\lambda\in\mathcal{C},
  \end{aligned}\right.
  \end{equation}
  where $\mathcal{C}$ represents the set of all complex numbers. Due to $\sqrt{Q_{k-1}}X = {\bf 0 }$ and $R>0$, one has $\sqrt{Q} X= {\bf 0}$ and $ K_{k-1}X={\bf0}$. As a result, equation (\ref{pf1eq1}) can be reduced to 
\begin{equation}
\left\{\begin{aligned}
&XA+AX+CXC^{\mathrm{T}} =\lambda X,\\
&\sqrt{Q}X={0},~~\lambda\in\mathcal{C},
\end{aligned}\right.
\end{equation}
which implies that $[A,C;\sqrt{Q}]$ is not exactly observable, contradicting the basic assumption ({\bf A3}).

By the obtained exact observability of $[A_{k-1},C_{k-1};\sqrt{Q_{k-1}}]$ and the induced assumption that $K_{k-1}$ is a stabilizer, according to \cite[Theorem 6]{zhang2004stabilizability}, equation (\ref{pk}) admits a unique solution $P_{k}>0$.

To establish that $K_{k}$ is a stabilizer, we use (\ref{pk}) and (\ref{kk}) to derive the following equation
\begin{equation}\label{pf0eq1}
\begin{aligned}
 A_{k}^{\mathrm{T}} P_{k} + P_{k} A_{k}+C_{k}^{\mathrm{T}} P_{k} C_{k}+\tilde{Q}_{k}={\bf0},
\end{aligned}
\end{equation}
where $\tilde{Q}_{k}= Q+(K_{k-1}-K_{k})^{\mathrm{T}}(R+D^{\mathrm{T}}P_{k}D)(K_{k-1}-K_{k})+K_{k}^{\mathrm{T}}RK_{k}$.

Similar to the above derivation, it is easy to verify that $[A_{k},C_{k};\sqrt{\tilde{Q}_{k}}]$ is exactly observable, combined with $P_{k}>0$, which results in $K_{k}$ being a stabilizer of $[A,B;C,D]$.

 Since $K_0$ is a stabilizer, the induction assumption is satisfied when $k=1$. Therefore, we establish results 1) and 2) for all $k\in\mathbb{N}_+$.

Next, we show that the sequence $\{P_{k}\}_{1}^{\infty}$ is monotonically decreasing with a lower bound $P$. 

To proceed, in view of (\ref{kk}), (\ref{pk}), and (\ref{ARE1}), we have
\begin{equation}\label{pf0eq2-1}
\begin{aligned}
 &A_{k}^{\!\mathrm{T}}(P_k\!-\!P_{k+1})\!+\!(P_{k}\!-\!P_{k+1}\!)A_{k}\!+\!C_{k}^{\!\mathrm{T}}(P_{k}\!-\!P_{k+1}\!)C_{k}\\
 =&\!-(K_{k-1}-K_{k})^{\mathrm{T}}\!(R+D^{\mathrm{T}}P_kD)\!(K_{k-1}-K_k),
\end{aligned}
\end{equation}
and
\begin{equation}\label{pf0eq2}
\begin{aligned}
&A_{\!k-1}^{\mathrm{T}}(P_k\!-\!P)\! + \!(P_k\!-\!P)A_{\!k-1}\!+\!C_{\!k-1}^{\!\mathrm{T}}(P_k\!-\!P)C_{k-1}\\
=&-(K_{k-1}-K)^{\mathrm{T}}\Upsilon(K_{k-1}-K).
\end{aligned}
\end{equation}
As $K_{k-1}$ is a stabilizer for all $k\in\mathbb{N}_+$, in conjunction with the conditions $(\!K_{k-1}\!-\!K)^{\!\mathrm{T}}\!(R+D^{\mathrm{T}}P_kD)\!(K_{k-1}\!-\!K)\!\geq\!{0}$ and $(\!K_{k-1}\!-\!K)^{\!\mathrm{T}}\!\Upsilon\!(K_{k-1}\!-\!K)\!\geq\!{0}$, equations (\ref{pf0eq2-1}) and (\ref{pf0eq2}) give rise to the following inequalities:
\begin{equation}
P_{k}\geq P_{k+1},P_k\geq P,~~\forall k\in\mathbb{N}_+,
\end{equation} 
 Therefore, $\{P_k\}_1^{\infty}$ is a monotonically decreasing sequence with a lower bound $P$. Since $(P,K)$ satisfies (\ref{pk}), it can be concluded that $\lim_{k\rightarrow\infty}P_k=P$.

Hence, the proof is complete.\hfill $\square$

% \end{pf}
     % \bibliographystyle{unsrt}
\bibliographystyle{bbs}        % Include this if you use bibtex 
\bibliography{autosam}                                    
\end{document}